\documentclass[11pt, letterpaper, oneside]{amsart}

\usepackage{latexsym, amsmath, amssymb, amsfonts, amscd}
\usepackage{amsthm}
\usepackage{t1enc}
\usepackage[mathscr]{eucal}
\usepackage{indentfirst}
\usepackage{graphicx, pb-diagram}
\usepackage{fancyhdr}
\usepackage{fancybox}
\usepackage{enumerate}
\usepackage[all]{xy}

\theoremstyle{plain}
\newtheorem{thm}{Theorem}[section]

\newtheorem{prop}[thm]{Proposition}
\newtheorem{lemma}[thm]{Lemma}
\newtheorem{cor}[thm]{Corollary}
\theoremstyle{definition}
\newtheorem{defi}[thm]{Definition}

\theoremstyle{remark}
\newtheorem{remark}[thm]{Remark}
\newtheorem{ep}[thm]{Example}

\newcommand{\ZZ}{\ensuremath{\mathbb Z}}

\newcommand{\RR}{\ensuremath{\mathbb R}}
\newcommand{\g}{\ensuremath{\frak{g}}}

\newcommand{\cP}{\mathcal{P}}
\newcommand{\cL}{\mathcal{L}}
\newcommand{\cE}{\mathcal{E}}
\newcommand{\cF}{\mathcal{F}}

\newcommand{\tn}{\tilde{\nabla}}
\newcommand{\tg}{\tilde{\gamma}}
\newcommand{\tv}{\tilde{v}}

\newcommand{\fG}{{f_\Gamma}}

\newcommand{\bt}{\mathbf{t}}                  %target
\newcommand{\bs}{\mathbf{s}}                  %source

\newcommand{\cG}{\mathcal{G}}
\newcommand{\Gs}{\Gamma_s}
\newcommand{\Gc}{\Gamma_c}

\newcommand{\bL}{\bar{L}}

\newcommand{\ta}{\tilde{a}}
\newcommand{\tP}{\tilde{P}}
\newcommand{\tQ}{\tilde{Q}}
\newcommand{\vG}{v_{\Gamma}}
\newcommand{\EG}{E_{\Gamma}}
\newcommand{\tG}{\theta_{\Gamma}}

\newcommand{\hEG}{\hat{E}_{\Gamma}}

\begin{document}

\title{On the geometry of prequantization spaces}
\author{Marco Zambon and Chenchang Zhu}
\date{\today}

\begin{abstract}
Given a Poisson (or more generally Dirac) manifold $P$, there are
two approaches to its geometric  quantization: one involves a
circle bundle $Q$ over $P$ endowed with a Jacobi (or Jacobi-Dirac)
structure; the other one involves a circle bundle with a 
(pre)contact groupoid structure over the (pre)symplectic groupoid of
$P$. We study the relation between these two prequantization
spaces. We show that the circle bundle over the (pre)symplectic
groupoid of $P$ is obtained from
 the Lie groupoid of $Q$ via an
$S^1$ reduction that preserves both the Lie groupoid and the geometric
structures.
\end{abstract}

 \maketitle \tableofcontents

\section{Introduction}

The geometric quantization of symplectic manifolds is a classical
problem that has been much studied over years. The first step is
to find a prequantization. A symplectic manifold $(P,\omega)$ is
prequantizable iff $[\omega]$ is an integer cohomology class.
Finding a prequantization means finding a faithful representation
of the Lie algebra of functions on $(P,\omega)$ (endowed with the
Poisson bracket) mapping the function $1$ to a multiple of the
identity. Such a representation space consists usually of sections
of a line bundle over $P$ \cite{kost}, or equivalently of
$S^1$-antiequivariant complex functions on the total space $Q$ of
the corresponding circle bundle \cite{So}. 

For more general kinds of geometric structure on $P$, such as
Poisson or even more generally  Dirac \cite{Co} structures, there
are two approaches to extend the geometric quantization of
symplectic manifolds, at least as far as prequantization is
concerned:
\begin{itemize}
\item To build a circle bundle $Q$ over $P$ compatible with the
Possion (resp. Dirac) structure on $P$ (see Souriau \cite{So}
for the symplectic case, \cite{Hu}\cite{va}\cite{cmdl} for the
Poisson case, and \cite{WZ} for the Dirac case)
\item To build the
symplectic (resp. presymplectic) groupoid of $P$ first and construct a
circle bundle over the groupoid \cite{wx}, with the hope to quantize
Poisson manifolds ``all at once'' as proposed by Weinstein
\cite{w-ncomquan}.
\end{itemize}

We call $Q$ as above a
``prequantization space'' for $P$ because, when $P$ is prequantizable, out of the hamiltonian vector fields on 
$Q$
one can construct  a representation of the admissible functions on $P$,
which form
 a Poisson algebra,
 on the space of $S^1$ anti-equivariant functions on $Q$ (see Prop. 5.1 of \cite{WZ}).
Usually however this representation is not faithful.

Since the (pre)symplectic groupoid $\Gamma_s(P)$ of $P$ is the canonical global object
associated to $P$, the prequantization circle bundle over $\Gamma_s(P)$ can 
 be considered an
``alternative prequantization space'' for $P$. Furthermore, since there is a
submersive
Poisson (Dirac) map $\Gamma_s(P)\rightarrow P$,  the admissible functions on $P$ can 
be viewed as a Poisson subalgebra of the functions on $\Gamma_s(P)$, which can be prequantized whenever $\Gamma_s(P)$ is a prequantizable (pre)symplectic manifold.
The resulting representation is faithful but the representation space is
unsuitable because  much too large.\\
 
In this paper we will not be interested in representations but only in the  
geometry that arises from the prequantization spaces associated to a  given a Dirac manifold $(P, L)$.
Indeed our main aim   is to study the relation between the  two
prequantization spaces above, which we will explain
in Thm. \ref{red}, Thm. \ref{redDirac} and Thm. \ref{gen-2}.

 We start searching for a more transparent description of the geometric
structures  on
the circle bundles $Q$, which are Jacobi-Dirac structures \cite{WZ} $\bar{L}$ .
This will be done in Section
\ref{cpreq}, both in terms of subbundles and in terms of brackets
of functions, paying particular attention  to the Lie algebroid structure that $\bar{L}$ carries.

Secondly, in Section \ref{recover}, we relate the Lie algebroid $\bar{L}$
associated to $Q$ to the Lie algebroid of the prequantization of $\Gamma_s(P)$. We do this using $S^1$ precontact reduction,
paralleling one of the motivating examples of symplectic
reduction: $T^*M//_0 G =T^*(M/G)$. This gives us  evidence at the
infinitesimal level for
 the relation between the Lie groupoid associated to  $Q$ and the prequantization of
 $\Gamma_s(P)$.
The latter relation between Lie groupoids will  be described in Section \ref{grpic},
again as an $S^1$ precontact
reduction.
 We provide a direct proof in the Poisson case. In the general Dirac case, the proof is
 done by integrating the results of Section \ref{recover} to the level of Lie groupoids with
 the help of Lie  algebroid path spaces. As a byproduct, we obtain the prequantization
condition for  $\Gamma_s(P)$
in terms of period groups on $P$. Then we show that this condition
is automatically satisfied when the Dirac manifold $P$ admits a  prequantization
  circle bundle $Q$ over it.
  This generalizes some of the results in \cite{cz}
and \cite{bcz}. 

This paper ends with three appendices. Appendix \ref{preciso} provides a
useful tool to perform computations on precontact groupoids, and
Appendix \ref{lcsgroid} describes explicitly the Lie groupoid of a locally
conformal symplectic manifold. In Appendix \ref{vor} we apply a
construction of Vorobjev to the setting of Section \ref{cpreq}.\\

\textbf{Notation:} Throughout the paper, unless otherwise
specified, $(P,L)$ will always denote a Dirac manifold,
$\pi:Q\rightarrow P$ will be a circle bundle and $\bar{L}$ will be
a Jacobi-Dirac structure on $Q$. By $\Gamma_s$ and $\Gamma_c$ we
will denote presymplectic and precontact groupoids respectively,
and we adopt the convention that the source map induces the (Dirac
and Jacobi-Dirac respectively) structures on the bases of the
groupoids. By ``precontact structure'' on a manifold we will just
mean a 1-form on the manifold.\\

\textbf{Acknowledgements:}  M.Z. is indebted to Rui Fernandes, for
an instructive invitation to IST Lisboa in January 2005, as well
as to Lisa Jeffrey. C.Z. thanks Philip Foth, Henrique Bursztyn and
Eckhard Meinrenken  for invitations to their institutions. Both
authors are indebted to Alan Weinstein for his invitation to U.C.
Berkeley in February/March 2005 and to the organizers of the
conference GAP3 in Perugia (July 2005).
 Further, we
thank A. Cattaneo and K. Mackenzie for helpful discussions, and
Rui Fernandes for suggesting the approach used in Subsection
\ref{fct} and pointing out the reference \cite{Vo}.

\section{Constructing the prequantization of $P$}\label{cpreq}

The aim of this section is to describe in an intrinsic way the
geometric structures (Jacobi-Dirac structures $\bar{L}$) on the
circle bundles $Q$ induced by prequantizable Dirac manifolds
$(P,L)$, paying particular attention to the associated Lie algebroid structures.
In Subsection \ref{alan} we will recall the non-intrinsic construction of $\bar{L}$
given in \cite{WZ}.
 In Subsection \ref{algoid} we will  describe $\bar{L}$ intrinsicly in
terms of subbundles and in Subsection \ref{fct}
  by specifying the bracket on functions
that it induces.\\

We first recall few definitions from \cite{WZ}. 
\begin{defi}
A \emph{Dirac structure}
on a manifold $P$ is  a subbundle of
$TP\oplus T^*P$ which is  maximal isotropic w.r.t. the symmetric pairing
$\langle X_1\oplus \xi_1,X_2\oplus \xi_2\rangle _+ = \frac{1}{2}
(i_{X_2}\xi_1+ i_{X_1}\xi_2)$ and whose sections are closed under 
   the Courant bracket
$$[X_1\oplus\xi_1,X_2\oplus\xi_2]_{Cou}=
\big([X_1,X_2]\;\oplus\;\cL_{X_1}\xi_2-\cL_{X_2}\xi_1+\frac{1}{2}d(
i_{X_2}\xi_1 -i_{X_1}\xi_2 )\big).$$
\end{defi}
 If $\omega$ is a 2-form on $P$ then its graph
$\{X\oplus \omega(X,\bullet):X\in TP\}$ is a Dirac structure iff $d\omega=0$.
  Given a bivector 
$\Lambda$ on $P$,  
the graph
$ \{ \Lambda(\bullet,\xi)\oplus \xi: \xi \in T^*P\}$ is a Dirac structure iff
  $\Lambda$ is a Poisson bivector. A Dirac structure $L$ on $P$ gives rise to (and is encoded by) a singular foliation of $P$, whose leaves are endowed with presymplectic forms.
  
  A function $f$ on a Dirac manifold $(Q,L)$ is \emph{admissible}
 if there
exists a smooth vector field $X_f$ such that $X_f\oplus df$ is a
section of $L$. A vector field $X_f$ as above is called a
\emph{hamiltonian vector field} of $f$. The set of admissible functions,
with the bracket  $\{f,g\}=X_g\cdot f$, forms a Lie (indeed a Poisson) algebra.
Given a map $\pi:Q \rightarrow P$ and a Dirac structure $L$ on $Q$, for
every $q\in Q$ 
one can define the subspace $(\pi_{\star}L)_{\pi(q)}:= 
\{\pi_*X\oplus \mu:X \oplus \pi^*\mu \in
L_q\}\}$
of $T_{\pi(q)}P\oplus T_{\pi(q)}^*P$. Whenever $\pi_{\star}L$ is a well-defined and  smooth subbundle of $TP\oplus T^*P$ it is automatically a Dirac structure on $P$. In this case $\pi:(Q,L) \rightarrow (P, \pi_{\star}L)$ is said to be a \emph{forward Dirac map}. Similarly, if $P$ is endowed with some Dirac structure $L$,
 $(\pi^{\star}L)(q):= 
\{Y \oplus \pi^* \xi:\pi_*Y\oplus \xi \in L_{\pi(q)}\}$ (when a smooth subbundle)
defines a Dirac structure on $Q$, and  $\pi:(Q,\pi^{\star}L) \rightarrow (P,  L)$ is said to be a \emph{backward Dirac map}.

 \begin{defi}
 A \emph{Jacobi-Dirac structure}  on $Q$ is defined as
  a subbundle   of  $\cE^1(Q):=(TQ\times \RR)\oplus (T^*Q\times \RR)$
  which is maximal isotropic w.r.t. the symmetric pairing
  $$\langle (X_1,f_1)\oplus(\xi_1,g_1)\;,\;(X_2,f_2)\oplus(\xi_2,g_2) \rangle_+
= \frac{1}{2}(i_{X_2}\xi_1
 +i_{X_1}\xi_2+g_2f_1+g_1f_2)$$  
  and whose space of sections is
closed under the extended Courant bracket on $\cE^1(Q)$ given by  
\begin{equation}\label{extcoubra}
\begin{split}
[(X_1,f_1)\oplus(\xi_1,g_1)\;&,\;(X_2,f_2)\oplus(\xi_2,g_2)]_{\cE^1(Q)}=
\big([X_1,X_2],X_1\cdot f_2-X_2\cdot f_1\big)\\
&\oplus \big(\cL_{X_1}\xi_2-\cL_{X_2}\xi_1
+\frac{1}{2}d(i_{X_2}\xi_1
-i_{X_1}\xi_2)\\
&+
f_1\xi_2-f_2\xi_1+\frac{1}{2}(g_2df_1-g_1df_2-f_1dg_2+f_2dg_1),\\
&X_1\cdot g_2-X_2\cdot
g_1+\frac{1}{2}(i_{X_2}\xi_1-i_{X_1}\xi_2-f_2g_1+f_1g_2)\big).
\end{split}\end{equation}
\end{defi}
We mention two examples.
Given any 1-form (precontact structure) $\sigma$ on $Q$,
$\text{Graph} \left(
\begin{smallmatrix} {d\sigma} & \sigma \\ -\sigma & 0
\end{smallmatrix} \right) \subset \cE^1(Q)$
is a Jacobi-Dirac structure.  
Given a bivector field
$\Lambda$ and a vector field $E$  on $Q$ and with the notation
  $\tilde{\Lambda}\xi:=\Lambda(\bullet,\xi)$,  $\text{Graph} \left(
\begin{smallmatrix} \tilde{\Lambda} & -E \\ E & 0
\end{smallmatrix} \right) \subset \cE^1(Q)$ is  a Jacobi-Dirac
structure iff $(\Lambda, E)$ is a 
  Jacobi structure, i.e. by definition if it   
satisfies the Schouten bracket
conditions $[E,\Lambda]=0$ and $[\Lambda,\Lambda]=2E\wedge
\Lambda$. Further to a  Dirac structure  $L\subset TQ\oplus T^*Q$  there is an associated
Jacobi-Dirac structure  $$L^c:=
\{(X,0)\oplus(\xi,g):(X,\xi)\in L, g\in \RR\}\subset \cE^1(Q).$$
   
A function $f$ on a Jacobi-Dirac manifold $(Q,\bar{L})$ is \emph{admissible} if there exists a smooth vector field $X_f$ and a
smooth function $\varphi_f$ such that
$(X_f,\varphi_f)\oplus(df,f)$ is a section of $\bar{L}$, and $X_f$ is called a \emph{hamiltonian
vector field} of $f$. The set of admissible functions,  denoted by
$C^{\infty}_{adm}(Q)$, together with the bracket
$\{f,g\}=X_g\cdot f+f\varphi_g$
forms a Lie algebra. There is a notion of forward and backward Jacobi-Dirac maps analogous to the one for Dirac structures.

\begin{defi}
A \emph{Lie algebroid} over a manifold $P$ is a vector bundle
$A$ over $P$ together with a Lie bracket
$[\cdot,\cdot]$ on its space of sections and a bundle map $\rho:
A\rightarrow TP$ (the \emph{anchor}) such that the Leibniz rule
$[s_1,fs_2]=\rho s_1(f)\cdot s_2+f\cdot [s_1,s_2]$ is satisfied
for all sections $s_1,s_2$ of $A$ and functions $f$ on $P$.
\end{defi}
One can think of Lie algebroids as generalizations of tangent bundles. To 
every Lie algebroid $A$ one associates cochains (the sections of the exterior algebra of $A^*$) and a certain differential $d_A$; the associated \emph{Lie algebroid cohomology} $H_A^{\bullet}(P)$ can be thought of as a generalization of deRham cohomology. One also defines an \emph{$A$-connection}  on a vector bundle $K\rightarrow P$ as map $\Gamma(A)\times \Gamma(K)\rightarrow \Gamma(K)$
satisfying the usual properties of a contravariant connection.

A Dirac structure $L\subset TP\oplus T^*P$ is automatically a Lie algebroid over $P$,
with bracket on sections of $L$ given by  the Courant bracket and anchor the   projection $\rho_{TP}:L\rightarrow TP$. Similarly, a Jacobi-Dirac structure $\bar{L}\subset \cE^1(Q)$, with the   extended Courant bracket and 
projection onto the first factor as anchor, is a Lie algebroid.\\

\subsection{A non-instrinsic description of $\bar{L}$}\label{alan}

We now recall the
prequantization construction of \cite{WZ},
which associates to a Dirac manifold a circle bundle $Q$ with a Jacobi-Dirac structure.

 Let $(P,L)$ be a Dirac structure. We saw above that  $L$ is a Lie algebroid
with the restricted Courant bracket and
 anchor $\rho_{TP}: L
\rightarrow TP$ (which is just the projection onto the tangent
component).  This anchor gives a Lie algebra homomorphism from
$\Gamma(L)$  to $\Gamma(TP)$ endowed with the Lie bracket of
vector fields.  The pullback by the anchor therefore induces a map
$\rho_{TP}^*: \Omega^{\bullet}_{dR}(P,\RR)\rightarrow
\Omega^{\bullet}_{L}(P)$, the sections of the exterior algebra of
$L^*$, which descends to a map from de Rham cohomology to the Lie
algebroid cohomology $H^{\bullet}_L(P)$ of $L$.
 There is a distinguished class in $H^{2}_{L}(P)$:
 on $TP\oplus T^*P$ there is an anti-symmetric pairing given by
\begin{eqnarray}
\langle X_1\oplus \xi_1,X_2\oplus \xi_2\rangle _- = \frac{1}{2}(
i_{X_2}\xi_1- i_{X_1}\xi_2).
\end{eqnarray}
Its restriction $\Upsilon$ to $L$ satisfies $d_L \Upsilon=0$. The
\emph{prequantization condition} (which for Poisson manifolds was first
formulated by Vaisman) is
\begin{eqnarray}\label{cond0} [\Upsilon] =
\rho_{TP}^*[\Omega]
\end{eqnarray}
for some integer deRham 2-class $[\Omega]$. 
$\eqref{cond0}$ can be equivalently phrased as
\begin{eqnarray} \label{cond1}\rho_{TP}^*\Omega=\Upsilon+d_L \beta,
\end{eqnarray}
where $\Omega$ is a closed integral 2-form and $\beta$ a 1-cochain
for the Lie algebroid $L$, i.e. a section of $L^*$. Let
$\pi:Q\rightarrow P$ be an $S^1$-bundle with connection form
$\sigma$ having curvature $\Omega$; denote by $E$ the
infinitesimal generator of the $S^1$-action. In Theorem 4.1 of
\cite{WZ} $Q$ was endowed with the following geometric structure,
described in terms of the triple $(Q,\sigma,\beta)$:

\begin{thm} \label{thmpreq}
The subbundle $\bar{L}$ of $\cE^1(Q)$ given by the direct sum of
$$\{(X^H+\langle X\oplus \xi, \beta \rangle
E,0)\oplus(\pi^*\xi,0):
 X\oplus \xi\in L\}$$ and the line bundles generated by
 $(-E,0)\oplus(0,1)$ and
 $(-A^H,1)\oplus(\sigma-\pi^*\alpha, 0)$
is a Jacobi-Dirac structure on $Q$. Here, $A\oplus \alpha$ is an
isotropic section of $TP\oplus T^*P$ satisfying $\beta=2\langle
A\oplus \alpha, \,\cdot\, \rangle _+|_L$. Such a section always
exists, and the subbundle above is independent of the choice of
$A\oplus \alpha$.
\end{thm}

 We call $(Q,\bar{L})$ a
``prequantization space'' for $(P,L)$ because the assignment
$g\mapsto \{\pi^*g,\bullet\}=-X_{\pi^*g}$ is a representation of
$C^{\infty}_{adm}(P)$ on the space of $S^1$ anti-equivariant functions on $Q$ \cite{WZ}. 

 Triples $(Q,\sigma,\beta)$ as
above define a hermitian $L$-connection with curvature $2\pi i
\Upsilon$ on the line bundle $K$ corresponding to $Q$, via the formula
\begin{eqnarray}\label{deco}D_{\bullet}=\nabla_{\rho_{TP} \bullet}-2\pi i \langle
\bullet,\beta \rangle
\end{eqnarray}
 where $\nabla$ is the covariant connection
corresponding to $\sigma$ (Lemma 6.2 in
\cite{WZ}). We have

\begin{prop}
For a prequantizable Dirac manifold $(P,L)$, the
Jacobi-Dirac structure $\bar{L}$ constructed in Thm. \ref{thmpreq}
on $Q$ is determined by a choice of  hermitian $L$-connection on $K$  with
curvature $2\pi i \Upsilon$. 
\end{prop}
\begin{proof}
We described above how the triples $(Q,\sigma,\beta)$ used to construct $\bar{L}$ give rise to hermitian $L$-connections with curvature $2\pi i \Upsilon$. 
Conversely, all hermitian $L$-connections with curvature $2\pi i \Upsilon$ arise
from triples  $(Q,\sigma,\beta)$ as above (Proposition 6.1 in
\cite{WZ}). A short computation shows that the triples that define
the same $L$-connection as $(Q,\sigma,\beta)$ are exactly those of
the form $(Q,\sigma+\pi^*\gamma,\beta+\rho_{TP}^*\gamma)$ for some
1-form $\gamma$ on $P$, and that these triples all define the same
Jacobi Dirac structure $\bar{L}$ (Lemma 4.1 in \cite{WZ}; see also
the last comment in Sect. 6.1 there).
\end{proof}

In the next two subsections we will construct $\bar{L}$ directly from the $L$-connection.  We end this subsection by commenting on how the various
Jacobi-Dirac structure $\bar{L}$ defined above are related.

\begin{remark}\label{howmany}
Two $L$-connections on $K$ are {\em gauge equivalent} if the differ by
$d_L\phi$ for some function $\phi:P\rightarrow S^1$.
Gauge-equivalent $L$-connections $D$ on $K$ with curvature $2 \pi
i \Upsilon$ give rise to isomorphic Jacobi-Dirac structures:
denoting by $\Phi$ the bundle automorphism of $Q$ given by
$q\mapsto q\cdot \pi^*\phi$, using the proof of Proposition 4.1 in
\cite{WZ} one can show that if $D_2=D_1-2 \pi  i d_L\phi$ then
$(\Phi_*,Id)\oplus ((\Phi^{-1})^*,Id)$ is an isomorphism from the
Jacobi-Dirac structure induced by $D_1$ to the one induced by
$D_2$. (Alternatively one can check directly that for the bracket of
functions, which by Remark \ref{Jacbr} determine the Jacobi-Dirac structures,
$\Phi^*\{\cdot,\cdot\}_{D_2}=\{\Phi^*\cdot,\Phi^*\cdot\}_{D_1}$.
 The gauge-equivalence classes of $L$-connections
with curvature $2 \pi i \Upsilon$ are a principal homogeneous
space for $H^1_L(P,U(1))$ (see the proof of Prop. 6.1 in
\cite{WZ}).
\end{remark}

\begin{remark}\label{mor}
It's easy to see that the prequantization space $Q$ of a
prequantizable Dirac manifold $(P,L)$ can be endowed with various
non-isomorphic Jacobi-Dirac structures $\bL$. Even more is true:
 $(Q,\bL_1)$ and $(Q,\bL_2)$ will usually not even be Morita
equivalent, for any reasonable notion of Morita equivalence of
Jacobi-Dirac manifold (or of their respective precontact
groupoids). Indeed for $P=\RR$ with the zero Poisson structure,
choosing $(Q,\sigma,\beta)=(S^1\times \RR,d\theta,x\partial_x)$ as
in Example \ref{1dim} one obtains a Jacobi structure on $Q$ with
three leaves, whereas choosing $(S^1\times \RR,d\theta,0)$ one
obtains a Jacobi structure with uncountably many leaves (namely
all $S^1\times \{q\}$). On the other hand, one of the general
properties of Morita equivalence is to induce a bijection on the
space of leaves.
\end{remark}

\subsection{An intrinsic characterization of  $\bar{L}$}\label{algoid}

In this subsection we fix an $L$-connection $D$ on  the line bundle $K\rightarrow P$ with
curvature $2\pi i \Upsilon$ and   construct the Lie algebroid
$\bL$ from $L$ and $D$ directly. (In Prop. \ref{algebroid-s1-red}
 we will perform the inverse construction, i.e. we will recover
$L$ from $\bar{L}$). An alternative approach that works in particular cases is presented in Appendix \ref{vor}.

We begin with a useful lemma concerning flat Lie algebroid connections
(compare also to Lemma 6.1 in \cite{WZ}).
\begin{lemma}\label{flat}
Let $E$ be any Lie algebroid over a manifold $M$, $K$ a line bundle
over $M$, and $D$ a Hermitian $E$-connection on $K$. Consider the
central extension $E\oplus_{\eta}\RR$, where $2 \pi i \eta$ equals
the curvature of $D$; then $\tilde{D}_{(Y,g)}=D_Y+2 \pi i g$
defines an
 $E\oplus_{\eta}\RR$-connection on $K$ which is moreover flat.
\end{lemma}
\begin{proof}
One checks easily that $\tilde{D}$ is indeed a Lie algebroid
connection.   Recall that the bracket on $E\oplus_{\eta}\RR$ is defined
as $[(a_1,f_1), (a_2,f_2)]_{E\oplus_{\eta}\RR}=
([e_1,e_2]_{E}, \rho(a_1)f_2-\rho(a_2)f_1+\eta(a_1,a_2))$, where $\rho$
is the anchor,
 and that the 
  curvature of $\tilde{D}$ is
\[ R_{\tilde{D}}(e_1, e_2)s = \tilde{D}_{e_1} \tilde{D}_{e_2} s -\tilde{D}_{e_2} \tilde{D}_{e_1} s - \tilde{D}_{[e_1,e_2]}  s \]
for elements $e_i$ of $E\oplus_{\eta}\RR$ and $s$ of $K$.
  The flatness of $\tilde{D}$ follows by a straightforward
calculation.
\end{proof}

We will use of this construction, which is just a way to make
explicit the structure of a transformation algebroid (see Remark
\ref{Kir} below).
% The assumptions on the group $G$ and on the
%connection $\tilde{D}$ on a vector bundle are just made to ensure that $\tilde{D}$
%comes from a connection on the corresponding principal $G$ bundle.
\begin{lemma}\label{pullback}
Let $A$ be any Lie algebroid over a manifold $P$, $\pi_Q:Q\rightarrow
P$ a principle
%$G$ bundle over $P$ for a subgroup $G\subset
$SO(n)$-bundle, $\pi_K:K\rightarrow P$ the vector bundle
associated to the standard representation of $SO(n)$ on $\RR^n$,
and $\tilde{D}$ a flat $A$-connection on $K$ preserving its
fiber-wise metric.
 The $A$-connection induces a bundle map $h_Q:\pi_Q^*A \rightarrow
TQ$ (the ``horizontal lift'')
 that can be used to extend, by the Leibniz rule, the obvious
bracket on $SO(n)$-invariant sections of $\pi_Q^*A$ to all
sections of $\pi_Q^*A$. The vector bundle $\pi_Q^*A$, with this
bracket and $h_Q$ as an anchor, is a Lie algebroid over $Q$.
\end{lemma}
\begin{proof}
We first recall some facts from Section 2.5 in \cite{Fe}. The
$A$-connection $\tilde{D}$ on the vector bundle $K$ defines a map
(the ``horizontal lift'')
 $h_K:\pi_K^*A
\rightarrow TK$ covering the anchor $A\rightarrow TP$ by taking
parallel translations of elements of $K$ along $A$-paths. See Section
\ref{dirid} for the definition of $A$-paths.
 Explicitly, fix an $A$-path $a(t)$ with base path $\gamma(t)$, a
point $x\in \pi_K^{-1}(\gamma(0))$
 and let
$\tilde{\gamma}(t)$ the unique path in $K$ (over $\gamma(t)$)
starting at $x$
 with
$\tilde{D}_{a(t)}  \tilde{\gamma}(t)=0$.
 We can always write $\tilde{D}=\nabla_{\rho \bullet} -\tilde{\beta}$ where
$\nabla$ is a metric $TP$-connection on $A$ and $\tilde{\beta}\in
\Gamma(A^*)\otimes \mathfrak{so}(K)$; then $\nabla_{\rho a(t)}
\tilde{\gamma}(t) =\langle \tilde{\beta},a(t) \rangle
\tilde{\gamma}(t)$. Since the left hand side is the projection of
the velocity of
 $\tilde{\gamma}(t)$ along the Ehresmann distribution $H$
corresponding to $\nabla$, we obtain $\frac{d}{dt}
\tilde{\gamma}(t)=  (\frac{d}{dt}{\gamma}(t))^H + \langle
\tilde{\beta},a(t) \rangle  \tilde{\gamma}(t),$ so that
\begin{eqnarray}\label{lifth}
h_K(a(0),x):=\frac{d}{dt}|_{t=0} \tilde{\gamma}(t) =\rho(a(0))^H+
\langle \tilde{\beta},a(0) \rangle x.
\end{eqnarray}
Of course $h_K$ does not depend on $\nabla$ or $\tilde{\beta}$
directly, but just on $\tilde{D}$. By our assumptions $h_K$ is
induced by a ``horizontal lift'' for the principle bundle $Q$,
i.e. by a $SO(n)$-equivariant map $h_Q:\pi_Q^*A \rightarrow TQ$
covering the anchor of $A$. Since our $A$-connection $\tilde{D}$
is flat,
 the map that associates to a section $s$ of $A$ the vector field
$h_Q(\pi_Q^*s)$ on $Q$ is a Lie algebra homomorphism.

On sections $\pi_Q^*s_1$,
 $\pi_Q^*s_2$ of $\pi^*_QA$ which are pullbacks of sections of $A$
we define the bracket to be $\pi_Q^*[s_1,s_2]$, and we extend it
to all sections of $\pi_Q^*A$ by using $h_Q$ as an anchor and
forcing the Leibniz rule. We have to show that the resulting
bracket satisfies the Jacobi identity. Given sections $s_i$ of $A$
and a function $f$ on $Q$ one can show that the Jacobiator
$[[\pi_Q^*s_1,f\cdot \pi_Q^*s_2],\pi_Q^*s_3]+c.p.=0$ by using the
facts that the bracket on sections of $A$ satisfies the Jacobi
identity and that the correspondence $\pi_Q^*s_i\mapsto
h_Q(\pi_Q^*s_i)$ is a Lie algebra homomorphism. Similarly, the
Jacobiator of arbitrary sections of $Q$ is also zero due to fact
that $h_Q$ actually induces a homomorphism on \emph{all} sections
of $\pi_Q^*A$.
\end{proof}

\begin{remark}\label{Kir}
Using $h_K$ instead of $h_Q$ in the construction of the previous
lemma leads to a Lie algebroid structure on $\pi^*_KA\rightarrow K$.
As Kirill Mackenzie pointed out to us, $\pi^*_KA$ is just the
transformation algebroid arising from the Lie algebroid action of $A$
on $K$ given by the flat connection $\tilde{D}$. Similarly, the Lie
algebroid structure on $\pi^*_QA$ we constructed in the lemma is
the transformation algebroid structure coming from $h_Q$, which is viewed here   as
a Lie algebroid action of $A$ on $Q$.
\end{remark}

Now we come back to our original setting, where we consider the Lie
algebroid $L$ over $P$ and a hermitian $L$-connection $D$ on the line
bundle $K$ over $P$. Consider $L^c$, the Jacobi-Dirac structure on $P$ naturally associated to $L$. There is a canonical  isomorphism 
$L^c \rightarrow L\oplus_{\Upsilon}\RR ,\;\;
(X,0)\oplus(\xi,g) \mapsto    (X,\xi,g)$
of Lie algebroids over $P$ \cite{cz}.
 Lemma \ref{flat} provides us with
 a flat $L\oplus_{\Upsilon}\RR$-connection $\tilde{D}$ on $K$ , and
by Lemma \ref{pullback} the pullback of $L\oplus_{\Upsilon}\RR$ to
$Q$ (the circle bundle associated to $K$) is endowed with a Lie algebroid structure. Using equation
\eqref{lifth} one sees that its anchor
$h_Q:\pi_Q^*(L\oplus_{\Upsilon}\RR)\rightarrow TQ$, at any point
of $Q$, is given by
 \begin{equation}\label{anchor} h_Q(X,\xi,g)=X^H +
(\langle X\oplus \xi, \beta \rangle  -g)E
\end{equation}
 (here we make immaterial choices to write $D$  as in equation \eqref{deco}
 and denote by $^H$ the horizontal lift w.r.t. $\ker
\sigma$). This formula for the anchor suggests how to identify 
$\pi_Q^*(L\oplus_{\Upsilon}\RR)$ with a subbundle of $\cE^1(Q)$:
we will show that
 the natural injection $$I:\pi^*_Q
(L\oplus_{\Upsilon}\RR) \rightarrow \bL\subset \cE^1(Q),\;\;\; 
I(X,\xi,g)=(h_Q(X,\xi,g),0)\oplus(\pi^*\xi,g)$$
is a
Lie algebroid morphism, whose image is a codimension one
subalgebroid of $\bar{L}$ which we denote by $\bar{L}_0$.
We regard $\bar{L}_0$ as a ``lift'' of $L$ (or rather $L^c$) obtained using the 
hermitian $L$-connection $D$.
Now we can describe the Jacobi-Dirac structure $\bar{L}$ prequantizing $L$ in invariant terms and characterize partially (see also Remark \ref{computing}) its Lie algebroid structure:

\begin{thm}\label{summ}
Assume that the Dirac manifold $(P,L)$ satisfies the
prequantization condition \eqref{cond0}. Fix the line bundle $K$
over $P$ associated with $[\Omega]$ and a Hermitian $L$-connection
$D$ on $K$ with curvature $2\pi i \Upsilon$. Denote as above by
$\bar{L}_0$
 the lift of  $L^c$ by the connection $D$. 
 Then $\bar{L}$, the subbundle defined in Thm. \ref{thmpreq},  is characterized as the 
 unique Jacobi-Dirac structure on $Q$ which contains $\bar{L}_0$
  and which is different from  $(\pi^{\star}L)^{c}$ (where $\pi^{\star}L$ denotes the pullback Dirac structure of $L$). Further $\bar{L}_0$
  is canonically isomorphic to   $\pi^*_Q
(L\oplus_{\Upsilon}\RR) $ as a Lie algebroid.
 \end{thm}
 
\begin{proof}
We first show that $I:\pi^*_Q
(L\oplus_{\Upsilon}\RR) \rightarrow \bL$ is indeed a Lie algebroid morphism. We compute for $S^1$ invariant sections
\begin{equation}\label{sections}\begin{split}
 [I(X_1,\xi_1,0),\;&I(X_2,\xi_2,0)]_{\cE^1(Q)}\\
 =&I([(X_1,\xi_1),(X_2,\xi_2)]_{Cou},0)+\langle (X_1,\xi_1),(X_2,\xi_2)
\rangle _{-} \left((-E,0)\oplus(0,1)\right)\\
=&I([(X_1,\xi_1,0),(X_2,\xi_2,0)]_{\pi_Q^*(L\oplus_{\Upsilon}\RR)})
\end{split}
\end{equation}
 and
$[I(X,\xi,0),I(0,0,1)]_{\cE^1(Q)}=0$; then one checks that $I$
respects the anchor maps of $\pi^*_Q (L\oplus_{\Upsilon}\RR)$ and
$\bL$.

To prove the above characterization of $\bar{L}$ we show that 
there are exactly two maximally isotropic subbundles of $\cE^1(Q)$ containing $\bar{L}_0$. Indeed, denoting by $(\bar{L}_0)^{\perp}$ the orthogonal of $\bar{L}_0$
w.r.t. the pairing $\langle \bullet,\bullet \rangle_+$, 
the quotient $(\bar{L}_0)^{\perp}/\bar{L}_0 $ is a rank 2 vector bundle over $Q$ which inherits  a non-degenerate symmetric pairing on its fibers. Every fiber of such bundle  is isomorphic to 
$\RR^2$ with pairing $\langle (a,b), (a',b') \rangle=\frac{1}{2} (ab'+ba')$,
which clearly contains exactly two isotropic subspaces of rank one (namely
$\RR(1,0)$ and $\RR(0,1)$). So there are at most two maximally isotropic subbundles of 
$\cE^1(Q)$ containing $\bar{L}_0$; indeed there are exactly two:
$\bar{L}$ and $\bar{L}_0\oplus \RR ((0,0)\oplus(0,1))$. The latter is 
$\pi^{\star}L=\{Y \oplus \pi^* \xi:\pi_*(Y)\oplus \xi \in L\}$ viewed as a Jacobi-Dirac structure on $Q$, hence we are done.
   \end{proof}

\begin{remark}\label{isoform}
Using the canonical identifications of Lie algebroids
  $L\oplus_{\Upsilon}\RR \cong L^c$ and 
 $\pi^*_Q (L\oplus_{\Upsilon}\RR)\cong\bar{L}_0$
  the 
  natural Lie algebroid morphism  
  $\pi^*_Q (L\oplus_{\Upsilon}\RR) \rightarrow  L\oplus_{\Upsilon}\RR$ is  
 \begin{equation}\label{morp}\Phi: \bar{L}_0\rightarrow L^c, (X,0)\oplus(\pi^*\xi,g)\mapsto (\pi_*X,0)\oplus
 (\xi,g). \end{equation}
 \end{remark}

\begin{remark}
The construction of Thm. \ref{summ} gives a quick way to
see that the subbundle $\bar{L}$ of $\cE^1(Q)$, as defined in Thm.
\ref{thmpreq}, is indeed closed under the extended Courant bracket:
$\bar{L}_0$ is closed since we realized it as a Lie algebroid, and  the sum with the span of 
the section $(-A^H,1)\oplus(\sigma-\pi^*\alpha, 0)$ is closed under the bracket
because $\langle [s_1,s_2]_{\cE^1(Q)}, s_3 \rangle_+$ (for $s_i$ sections of $\cE^1(Q)$)
is a totally skew-symmetric tensor 
 \cite{IMa}.
\end{remark}

\begin{remark}\label{computing}
The characterization of $\bar{L}_0$ as the transformation algebroid of some action of $L\oplus_{\Upsilon}\RR \cong L^c$ on $Q$  (Thm. \ref{summ}) shows that 
if the Lie algebroid $L^c$ is integrable then $\bar{L}_0$  is integrated by the
corresponding transformation groupoid.  
Unfortunately using Thm. \ref{summ} we are not able to make the same conclusion for
 $\bL$.  Looking at the brackets on $\bar{L}$ is not very illuminating: 
it is determined by  \eqref{sections}  and
\begin{equation}\begin{split}
[I(X,\xi,0),&\;(-A^H,1)\oplus(\sigma-\pi^*\alpha,0)]_{\cE^1(Q)}=
I(-[(X,\xi),(A,\alpha)]_{Cou},0)\\
+& I(0,\Omega(X)-\xi+\frac{1}{2} d\langle X\oplus \xi, \beta
\rangle,0 )-\langle A ,\xi\rangle \left((-E,0)\oplus(0,1)\right).
\end{split} \end{equation}
The remaining brackets between sections of the form $I(X,\xi,0)$,
$I(0,0,1)$ and $(-A^H,1)\oplus(\sigma-\pi^*\alpha,0)$ vanish, and
by the Leibniz rule these brackets determine the bracket for
arbitrary sections of $\bar{L}$.  \end{remark}

\begin{remark}
Different choices of $L$-connection on the line bundle $K$ with
curvature $2 \pi i\Upsilon$ usually lead to Lie algebroids $\bL$
with different foliations (see Remark \ref{mor}), which therefore
can not be isomorphic. However the subalgebroids $\bar{L}_0$ are
always isomorphic. Indeed any two  connections with the same curvature are of the
form $D$ and $D'=D+2\pi i \gamma$, where $\gamma$ is a closed
section of $L^*$ (see Prop. 6.1 in \cite{WZ}). A computation using
$d_L\gamma=0$ shows that $(X,\xi)\oplus g \mapsto (X,\xi)\oplus
(g-\langle (X,\xi),\gamma \rangle)$ is a Lie algebroid
automorphism of $L\oplus_{\Upsilon}\RR$. Further this automorphism
intertwines the Lie algebroid actions \eqref{anchor} of
$L\oplus_{\Upsilon}\RR$ on $Q$ given by the ``horizontal lifts''
of the flat connections $\tilde{D}$ and $\tilde{D'}$. Hence the
transformation algebroids of the two actions are isomorphic, as is
clear from the description of Lemma \ref{pullback}.

 We exemplify the fact that actions coming from different flat connections are
 intertwined by a Lie algebroid automorphism (something that can
not occur if the anchor of the Lie algebroid is injective)
 in the simple case when the Dirac structure on $P$ comes from a
close 2-form $\omega$: the Lie algebroid action of
$TP\oplus_{\omega}\RR$ on $Q$ via a connection $\nabla$ (with
curvature $2\pi i \omega$) is intertwined to the obvious action of
the Atiyah algebroid $TQ/S^1$ on $Q$ (essentially given by the
identity map) via $TP\oplus_{\omega}\RR\cong TQ/S^1$ is
$(X,g)\mapsto X^H-\pi^*gE$, where $\sigma$ is the connection on
the circle bundle $Q$ corresponding to $\nabla$.
\end{remark}

\subsection{Describing $\bar{L}$ via the bracket on
functions}\label{fct}

In this subsection we will  describe the geometric
structure $\bar{L}$ on the circle bundle $Q$ in terms of the bracket on the
admissible functions on $Q$; by Remark \ref{Jacbr} below the bracket on 
functions   uniquely determines $\bar{L}$.

We adopt the following notation.  $F_S$ denotes the
 function on $Q$ associated to a section $S$ of the line bundle $K$:  $F_S$ is just the
restriction to the bundle of unit vectors $Q$ of the fiberwise
linear function on $K$ given by $\langle \cdot, S \rangle $, where
$\langle \cdot, \cdot \rangle $ is the
 $S^1$-invariant
\emph{real} inner product on $K$ corresponding to the chosen
Hermitian form on $K$.  Alternatively $F_S$ can be described as the
real part of the $S^1$-antiequivariant function on $Q$ that
naturally corresponds to the section $S$. By
 $iS$ we denote the
image of the section $S$ by the action of $i\in S^1$ (i.e. $S$
rotated by $90^{\circ}$), and $f$ and $g$ are functions on $P$.

%The following theo says how to endow circle bundles over
%$(P,L)$ with the geometric structure of a ``prequantization
%space''.
\begin{prop}\label{bracket}
Assume that the Dirac manifold $(P,L)$ satisfies the
prequantization condition \eqref{cond0}. Fix the line bundle $K$
over $P$ associated with $[\Omega]$ and a Hermitian $L$-connection
$D$ on $K$ with curvature $2\pi i \Upsilon$. Denote by $\tilde{D}$
the flat connection induced as in Lemma \ref{flat} and by
$h_Q:\pi_Q^*(L\oplus_{\Upsilon}\RR)\rightarrow TQ$ the horizontal
lift associated to $\tilde{D}$ given by \eqref{anchor}.

Suppose a Jacobi-Dirac structure $\hat{L}$ on $Q$ has the
following two properties: first, nearby any $q\in Q$ such that
$TP\cap L$ is regular near $\pi(q)$, the admissible functions for
$\hat{L}$ are exactly those that are constant along the leaves of
$\{h_Q(X,0,0):X\in TP\cap L\}$. Second, the bracket on locally
defined admissible functions is given by
\begin{itemize}
\item $\{\pi^*f,\pi^*g\}_Q=\pi^*\{f,g\}_P$
\item $\{\pi^*f,F_S\}_Q=F_{-\tilde{D}_{X_f,df,f}S}$
\item $\{\pi^*f,1\}_Q=0$
\item $\{F_S,1\}_Q=-2\pi F_{iS}$.
\end{itemize}
Then $\hat{L}$ must be the Jacobi-Dirac structure $\bar{L}$ given
in Thm. \ref{thmpreq}.

Conversely, the Jacobi-Dirac structure $\bar{L}$ given in Thm.
\ref{thmpreq} has the two properties above.
\end{prop}

\begin{proof}
We start by showing that the Jacobi-Dirac structure $\bar{L}$
constructed in Thm. \ref{thmpreq} satisfies the above two
properties. On the set of points where the ``characteristic
distribution'' $C:=\bar{L}\cap (TQ\times \RR)\oplus(0,0)$ of any
Jacobi-Dirac structure has constant rank the admissible functions
are exactly the functions $f$ such that $(df,f)$ annihilate $C$.
In our case $C=\{X^H+\langle \alpha,X \rangle E: X\in L\cap
TP\}=\{h_Q(X,0,0):X\in TP\cap L\}$ is actually contained in $TQ$,
so the admissible functions are those constant on the leaves of
$C$ as claimed.

Now we check that the four formulae for the bracket hold. The
first equation follows from the fact that the pushforward of
$\bar{L}$ is the Jacobi-Dirac structure associated to $L$ (see
Section 5 in \cite{WZ}).

For the second equation we make use of the formulae
$$E(F_S)=-2 \pi F_{iS}\;\;\;\;\;\;\text{ and } X^H(F_S)=F_{\nabla_X
  S},$$
where we make some choice to express $D$ as in equation
\eqref{deco} and $X^H$ denotes to horizontal lift of $X\in TP$
using the
  connection on $Q$ corresponding to the covariant derivative $\nabla$
  on $K$. Using these formulae we see
\begin{eqnarray*} \{\pi^*f,F_S\}_Q&=& - \langle dF_S, X^H_f+\langle  (X_f,df),\beta
\rangle  E -f E \rangle  \\
&=& F_{-\nabla_{X_f}S+2\pi i (\langle  (X_f,df),\beta
\rangle  -f}\\
&=&F_{-\tilde{D}_{X_f,df,f}S}.
\end{eqnarray*}
For the last two equations just notice that, since $(-E,0)\oplus
(0,1)$ is a section of $\bL$,
 the bracket of any admissible function with the constant function $1$
amounts to applying $-E$ to that function.\\

Now we show that if a Jacobi-Dirac structure $\hat{L}$ satisfies
the two properties in the statement of the proposition, then it
must be $\bar{L}$. By Remark \ref{Jacbr}, the bracket of
$dimQ-rkC+1$ independent functions
 at regular points of $C:=\hat{L}\cap
(TQ\times \RR)\oplus(0,0)$ determines $\hat{L}$, so we  have to
show that our two properties carry the information of the bracket
of $dimQ-rkC+1$ independent functions
 at regular points of $C$ .

It will be enough to consider the open dense subset of the regular
points of $C$ where $C=\{h_Q(X,0,0):X\in TP\cap L\}$ (This subset
is dense because it includes the points $q$ such that $C$ is
regular near $q$ and $TP\cap L$ is regular near $\pi(q)$). Since
there $C$ is actually contained in $TQ$
 it is clear that $1$ and $\pi^*f$ are admissible functions,
for $f$ any admissible function on $P$ (this means that $f$ is
constant along the leaves of $L\cap TP$; there are $dim P-rk C$
such $f$ which are linearly independent at $\pi(q)$). Further we
can construct an admissible function $F_S$ as follows: take a
submanifold $Y$ near $\pi(q)$ which is transverse to the foliation
given by $L\cap TP$, and define the section $S|_Y$ so that it has
norm one (i.e. its image lies in $Q\subset K$). Then extend $S$ to
a neighborhood of $\pi(q)$ by starting at a point $y$ of $Y$ and
``following'' the leaf of $C$ through $S(y)$ (notice that $C$ is a
flat partial connection on $Q\rightarrow P$ covering the
distribution $L\cap TP$ on $P$). Since $C$ is $S^1$ invariant, the
resulting function $F_S$ is clearly constant along the leaves of
$C$, hence admissible. Altogether we obtain $dim Q-rk C+1$
admissible functions in a neighborhood of $q$ for which we know
the brackets, so we are done.
\end{proof}

\begin{remark}\label{Jacbr}
On any Jacobi-Dirac manifold $(Q,\hat{L})$ the bracket on the
sheaf of admissible functions
$(C_{adm}^{\infty}(Q),\{\cdot,\cdot\})$ determines the subbundle
$\hat{L}$ of $\cE^1(Q)$. (This might seem a bit surprising at
first, since the set of admissible functions is usually much
smaller than $C^{\infty}(Q)$).

The set of points where $C:=\hat{L}\cap (TQ\times \RR)\oplus(0,0)$
(an analog of a ``characteristic distribution'') has locally
constant rank is an open dense subset of $Q$, since $C$ is an
intersection of subbundles. Hence by continuity it is enough to
reconstruct the subbundle $\bL$ on each point $q$ of this open
dense set.

Since we assume that $C$ has constant rank near $q$, given
$C_{adm}^{\infty}(Q)$ in a neighborhood of $q$ we can reconstruct
$C$ as the distribution annihilated by $(df,f)$ where $f$ ranges
over $C_{adm}^{\infty}(Q)$. We can clearly find $dimQ-rkC+1$
admissible functions $f_i$ such that $\{(df_i,f_i)\}$ forms a
basis of $\rho_{T^*Q\times
  \RR}(\hat{L})=C^{\circ}$ near $q$. The fact that each $f_i$ is
  an admissible function means that there exist  $(X_i,\phi_i)$
  such that  $(X_i,\phi_i)\oplus (df_i,f_i)$ is a smooth section
of $\hat{L}$. Now knowing the bracket of any $f_j$ with the other
$f_i$'s,  i.e. the pairing of $(X_j, \phi_j)$ with all elements of
$\rho_{T^*Q\times
  \RR}(\hat{L})$, does not quite determine $(X_j, \phi_j)$. However it determines
  $(X_j, \phi_j)$ up to sections of
  $C$, hence the direct sum
of the span of all $(X_i,\phi_i)\oplus (df_i,f_i)$ and of $C$ is a
well  defined  subbundle of $\cE^1(Q)$. Moreover it has the same
dimension as $\hat{L}$ and it is spanned by sections of $\hat{L}$,
so it is $\hat{L}$.
\end{remark}

\section{Prequantization and reduction of Jacobi-Dirac structures}\label{recover}

In the last section we considered a  prequantizable Dirac manifold
$(P,L)$ and endowed $Q$ (the total space of the circle bundle over
$P$) with distinguished Jacobi-Dirac structures $\bar{L}$.
  
We are interested in the relation between the \emph{Lie algebroid} structures on $\bar{L}$ 
and   
$L^c$ (the Jacobi-Dirac structure canonically associated to $L$), because they will give an indication of the relation between the Lie groupoids integrating them.
The map $\Phi$ of \eqref{morp} is a natural surjective morphism of Lie algebroids from 
  the codimension one subalgebroid $\bar{L}_0$ of $\bar{L}$
to  $L^c$, so one may hope to extend  
 $\Phi$   to a Lie algebroid
morphism defined on $\bar{L}$. However in general there cannot be any
Lie algebroid morphism from $\bar{L}$ to $L^c$ or $L$ with base
map $\pi$: recall that a morphism of Lie algebroids maps each orbit of
the source Lie algebroid into an orbit of the target Lie algebroid. If the
map $\pi:Q\rightarrow P$ induced a morphism of Lie algebroids, then
the orbits\footnote{The orbits of a Lie algebroid are
 the leaves integrating the (singular)
 distribution given by
    the image of the anchor map.} of $\bar{L}$ would be mapped into the orbits of $L^c$
(which coincide with those of $L$). However this happens exactly
when (one and hence all choices of) the vector field $A$ appearing
in Thm. \ref{thmpreq} is tangent to the foliation of $L$ (see
Section 4.1 of \cite{WZ}). In the case of
 Example \ref{1dim}, i.e. $Q=S^1\times \RR $ and $P=\RR$, the orbits of $T^*Q\times
\RR$  are exactly three (namely $S^1\times \RR_+, S^1\times \{0\}$
and $S^1\times \RR_-$), and $\pi$ does not map them into the
orbits of $T^*P$, which are just points.
 
 In this section we will take advantage of the fact that $\bar{L}$, in addition to the Lie algebroid structure, also carries a geometric structure,
 namely a precontact structure $\theta_{\bL}\in \Omega^1(\bar{L})$ defined as follows:
\begin{equation}\label{precont-form}
\theta_{\bL}:=pr^* (\theta_c+dt),
\end{equation}
where $\theta_c$ is the canonical 1-form on the cotangent bundle $T^*Q$, $t$ is the coordinate on $\RR$, and
  $pr$ is the projection of
$\bL\subset \cE^1(Q)$ onto $T^*Q\times \RR$.
We will use the
the 1-form $\theta_{\bar{L}}$
to
recover the Lie
algebroid $L^c$ from $\bar{L}$ via a precontact reduction procedure, which we
will globalize to the corresponding Lie groupoids in the next
Section.

\subsection{Reduction of Jacobi-Dirac structures as precontact reduction}

We recall a familiar fact: in symplectic geometry, we have the
well-known motivating example of symplectic reduction $T^*M//_0 G=
T^*(M/G)$. In \cite{tudor}, it is extended to contact geometry by
replacing $T^*M$ by the cosphere bundle of $M$. Here we prove a
similar result by replacing $T^*M$ by $T^*M \times \RR$--another
natural contact manifold associated to any manifold $M$. Later on we will use this to reduce a $G$-invariant Jacobi-Dirac structure on $M$ to a Jacobi-Dirac structure on $M/G$.\\

 Let a Lie
group $G$ act on a contact manifold $(C, \theta)$ preserving the
contact form $\theta$. 

 Then, a moment map  is a map $J$ from the
manifold $M$ to $\g^*$ (the dual of the Lie algebra)  such that
for all $v$ in the Lie algebra $\g$:
\begin{equation}\label{cmmap} %contact moment map
 \langle J, v \rangle = \theta_M (v_M),
\end{equation}
where $v_M$ is the infinitesimal generator of the action on $M$
given by $v$. The moment map $J$ is automatically equivariant with
respect to the coadjoint action of $G$ on $\g^*$ given by $\xi
\cdot g=L_g^* R_{g^{-1}}^* \xi$. A group action as above together
with its moment map is called {\em Hamiltonian}. Notice that any
group action preserving the contact form is Hamiltonian. In the above setting
there are
two ways to perform contact reduction, developed  by Albert \cite{albert} and Willett 
  \cite{willett} respectively, which agree when one performs reduction at $0\in \g^*$:  \[ C//_0 G := J^{-1}(0)/ G \]
is again a smooth contact manifold with   induced contact form
$\bar{\theta}$ such that
$\pi^*(\bar{\theta})=\theta|_{J^{-1}(0)}$.

\begin{lemma} \label{t*m} Let the group $G$ act on manifold $M$ freely and properly. Then $G$ has an induced action on the contact manifold
$(C:=T^*M\times \RR, \theta:=\theta_c + dt)$ where $\theta_c$ is
the canonical 1-form on $T^*M$ and $t$ is the coordinate on $\RR$.
Then this action is Hamiltonian and the contact reduction at $0$ is
\[ T^*M\times \RR //_0 G = T^*(M/G)\times \RR.\]\end{lemma}
\begin{proof}
The induced $G$ action on $T^*M\times \RR$ is by $g\cdot (\xi, t)
= ((g^{-1})^* \xi, t)$, and it preserves the 1-form $\theta_c +
dt$.
 The projection
of this action on $M$ is the $G$ action on $M$ so it is also free
and proper. Then the moment map $J$ is determined by
\[ \langle J(\xi, t), v\rangle = (\theta_c +dt)_{(\xi,
t)}(v_C)=\theta_c(v_C)=\langle \xi, v_M\rangle  ,\] where $v_C$
(resp. $v_M$) denotes the vector filed corresponding to the
infinitesimal action of $G$ on the manifold $C$ (resp. $M$). Since
all infinitesimal generators $v_C$ are nowhere proportional to the
Reeb vector field $\frac{\partial}{\partial t}$,  by Remark 3.2 in
\cite{willett}  all points of $T^*M\times \RR$ are regular points
of $J$.
 So $J^{-1}(0)=\{ (\xi, t): \langle \xi,
v_M\rangle =0 \; \forall v\in \g\}=
\{(\pi^*\mu,t): \mu \in T^*(M/G)\}$ (with $\pi: M\to (M/G)$) is a smooth manifold.
Therefore it is not hard to see that there is a well-defined
\[\Phi: J^{-1}(0)/G\to T^*(M/G)\times \RR, \quad \; ([\xi],t) \mapsto
(\mu, t), \] where $\mu$ is uniquely determined by $\pi^* \mu=\xi$
and we used the notation $[\cdot]$ to denote the quotient of
points (and later tangent vectors) of $J^{-1}(0)$ by the $G$
action. It is not hard to see that $\Phi$ is an isomorphism since
the two sides have the same dimension and $\Phi$ is obviously
surjective. The contact form  on $T^*(M/G)\times \RR$
corresponding to the reduced contact form $\bar{\theta}$ via the
isomorphism $\Phi$
 is the canonical one: for a tangent vector $([v],\lambda \frac{\partial}{\partial t}
)\in T_{[\xi], t}(J^{-1}(0)/G)$,
%=T_{[\xi], t} (T^*(M/G)\times \RR)$
%, where $[X]$ is the quotient
%of a tangent vector $X\in T_\xi (T^*M)$,
\[ \bar{\theta}_{[\xi], t}( [v], \lambda \frac{\partial}{\partial
t})=\theta_{\xi, t} (v, \lambda \frac{\partial}{\partial t}) =
\xi(p_*v) +\lambda=\mu(\bar{p}_*\Phi_*[v]) +\lambda,\]where $p:
T^*M \to M$ and $\bar{p}: T^*(M/G) \to M/G$. Here we used
$\bar{p}_*\Phi_*[v]=\pi_*p_* v$, which follows from the fact that
$\Phi$ is a vector bundle map, and we abuse notation by denoting
with the same symbol  a restriction of $\Phi$.
\end{proof}

This result extends to the precontact situation:
instead of the contact manifold $T^*M \times \RR$  we consider a Jacobi-Dirac subbundle
  $\bL\subset \cE^1(M)$, which together with the 
  1-form
  $\theta_{\bL}\in \Omega^1(\bar{L})$ defined in
\eqref{precont-form} is a precontact manifold.

\begin{prop}\label{algebroid-red}
When $(Q, \bL)$ is a Jacobi-Dirac manifold, $\bL$ is a precontact
manifold as described above. If the  group $G$ acts freely and
properly  on $Q$ preserving the Jacobi-Dirac structure, the action
lifts to a free proper Hamiltonian action on $\bL$ with moment map
$J$,
\[ \langle J((X,f)\oplus (\xi,g)), v\rangle =\theta_{\bL}|_{(X,f)\oplus(\xi,g)} (v_{\bL})
=\xi(v_Q). \] Write $\g_Q$ as a short form for $\{v_Q:v\in
\g\}\subset TQ$, and let $\pi_{\star} \bar{L}
\subset \cE^1(P)$ be the pushforward
of $\bar{L}$ via $\pi:Q\rightarrow P:=Q/G$.
 Then
\begin{enumerate}
\item   $J^{-1}(0)$ is a
subalgebroid of $\bL$ iff $\bL \cap (\g_Q,0)\oplus(0,0)$ has constant
rank, and in that case  $\bL//_0 G:=
J^{-1}(0)/G$  has an induced Lie algebroid structure;
%\item   the  pushforward $L_P$ of
%$\bL$ via $\pi:Q\rightarrow P$ is an induced Jacobi-Dirac
%structure on $P$, and the precontact reduction $\bL//_0 G:=
%J^{-1}(0)/G$ is canonically isomorphic to $L_P$ as Lie algebroid,
%iff $\bL \cap (\g_Q,0)\oplus(0,0)=\{0\}$;
\item  $J^{-1}(0)/G\cong \pi_{\star} \bar{L}$
 both as Lie algebroids and precontact manifolds,
iff $\bL \cap (\g_Q,0)\oplus(0,0)=\{0\}$. Here the precontact forms are the reduced 1-form on $J^{-1}(0)/G$
and the one defined as in \eqref{precont-form} on $\pi_{\star} \bar{L}$ respectively.
\end{enumerate}
\end{prop}
\begin{proof} The $G$ action on $Q$ lifts to $\bL$ by
$g\cdot ( X, f)\oplus ( \xi , g) = (g_* X , f)\oplus ( (g^{-1})^*
\xi, g)$, and the resulting moment map $J$ is clearly as claimed
in the statement.

 To prove $(1)$ we start with some linear algebra and fix
$x\in Q$. We have a map $\pi_*: T_xQ\rightarrow T_{\pi(x)}(Q/G)$,
hence we can push forward $\bL_x$ to
$$(\pi_{\star} \bar{L})_{\pi(x)}=\{(\pi_*X,f)\oplus(\mu,g):(X,f)\oplus(\pi^*\mu,g)\in
\bL_x\}$$ to obtain a linear Jacobi-Dirac subspace of
$\cE^1(Q/G)_{\pi(x)}$. Since $\bL$ is $G$ invariant, doing this
at every $x\in Q$ we obtain a well defined subbundle of
$\cE^1(Q/G)$,
 which however
 %is not a Jacobi-Dirac structure because
 %it
  might fail to be smooth\footnote{For example it is not
smooth when $G=\RR$, $Q=\RR^2$, $v_Q= \frac{\partial}{\partial x}$
and $\bL$ is the graph of the 1-form $\frac{y^2}{2}dx$.}. We have
a surjective map
\begin{equation}\begin{split} \label{surj}
 & \Phi:J^{-1}(0)=\{(X,f)\oplus(\xi,g)\in
\bL:\xi=\pi^*\mu \text{ for some } \mu\in
T_{\pi(x)}(Q/G) \} \rightarrow \pi_{\star} \bar{L}\\
  & (X,f)\oplus (\xi,g)\mapsto
(\pi_*X, f)\oplus (\mu, g) \end{split}\end{equation}
 whose kernel is exactly
$J^{-1}(0)\cap(\g_Q,0)\oplus(0,0)$ (Notice that the map is well
defined for $\pi$ is a submersion). So $J^{-1}(0)$ has constant
rank iff $J^{-1}(0)\cap(\g_Q,0)\oplus(0,0)=\bL\cap
(\g_Q,0)\oplus(0,0)$ does. In this case
 it is easy to see that
$J^{-1}(0)$ is closed under the Courant bracket: the Courant
bracket of two sections of $J^{-1}(0)$ lie in $\bL$ (because $\bL$
is closed under the bracket), therefore one just  has to show that
its cotangent component is annihilated by $\g_Q$. By a
straight-forward computation this is true for $G$-invariant
sections, and by the Leibniz rule it follows for all sections of
$J^{-1}(0)$, i.e. $J^{-1}(0)$ is a subalgebroid. Clearly
$J^{-1}(0)/G$ becomes a Lie algebroid with the bracket induced from
the one on $J^{-1}(0)$ and anchor $([X],f)\oplus([\xi],g)\mapsto
\pi_*X$ (where $[\cdot]$ denotes the equivalence relation given by
the $G$ action).

To prove $(2)$ consider the map $\Phi$ above.  It induces an isomorphism of
 vector bundles over $P$ between $J^{-1}(0)/G$ and $\pi_{\star} \bar{L}$ iff
it is fiberwise injective, i.e. iff
 $\bL
\cap(\g_Q,0)\oplus(0,0)=\{0\}$. Since $J^{-1}(0)/G$ (being a
precontact reduction) is a smooth manifold and $J^{-1}(0)/G\cong
\pi_{\star} \bar{L}$ is point-wise a subbundle of $\cE^1(P)$, it follows that
$\pi_{\star} \bar{L}$ is a smooth vector bundle over $P$.
We are left with showing that $\Phi$ induces an isomorphism of Lie
 algebroids and precontact manifolds.
Using the fact
that operations appearing in
 the definition of Courant bracket such as taking Lie derivatives commute with taking quotient
of $G$ (for example $\pi^*(L_{\pi_*X} \mu)=L_{X}\pi^*\mu$)
 we deduce that
$\Phi:J^{-1}(0) \rightarrow \pi_{\star} \bar{L}$ is a surjective morphism of Lie
algebroids, hence the induced map $\Phi: J^{-1}(0)/G\rightarrow
\pi_{\star} \bar{L}$ an isomorphism of Lie algebroids.

The isomorphism of precontact manifolds follows from an entirely
similar argument as in Lemma \ref{t*m}. We consider a tangent
vector $([w],\kappa \frac{\partial}{\partial s})\oplus ([v],
\lambda \frac{\partial}{\partial t} )\in T_{([X], f)\oplus ([\xi],
g)}(J^{-1}(0)/G)$, then $\Phi(([X], f)\oplus ([\xi],
g))=(\pi_* X, f)\oplus (\mu, g)$,
% where $[v]$ is the
%quotient of a tangent vector $v\in T_{X,f ,\xi, g} \bL$ and
where $\pi^*\mu=\xi$. So the induced 1-form $\bar{\theta}$ on $J^{-1}(0)/G$
satisfies,
\[ \bar{\theta}_{[X], f, [\xi], g}([w],\kappa \frac{\partial}{\partial
s})\oplus
 ([v], \lambda \frac{\partial}{\partial
t})=\theta_{X,f ,\xi, g} (w,\kappa \frac{\partial}{\partial
t})\oplus (v, \lambda \frac{\partial}{\partial t}) = \xi(p_* v)
+\lambda=\mu(\bar{p}_*\Phi_* [v]) +\lambda,\]where $p: \bL \to Q$
and $\bar{p}: \pi_{\star} \bar{L} \to P$ are projections.
Therefore
$\bar{\theta}=\Phi^*\theta_{\pi_{\star} \bar{L}}$ with $\theta_{\pi_{\star} \bar{L}}$ the canonical 1-form as in \eqref{precont-form}.
\end{proof}

\begin{remark} A special case of Prop.
 \ref{algebroid-red} is the usual reduction of basic 1-forms:
if the Jacobi-Dirac structure $\bar{L}$ of Prop.
 \ref{algebroid-red} comes from  1-form $\sigma$ on $Q$
such that $\g_Q\subset ker \sigma$, then the pushforward $\pi_{\star} \bar{L}$ is given by the unique 1-form $\sigma_{red}$ on $P=Q/G$ satisfying $\pi^* \sigma_{red}=\sigma$.
\end{remark}

\subsection{Reduction of prequantizing Jacobi-Dirac structures}\label{algpic}

Now we adapt the general theory of reduction of Jacobi-Dirac
manifolds  discussed  in the previous subsection to our situation, namely we
consider  a prequantization $Q$ of Dirac manifold $(P,L)$. Then
$Q$ is Jacobi-Dirac with a free and proper $S^1$ action which
preserves the Jacobi-Dirac structure $\bar{L}$. Let $L^c=
\{(X,0)\oplus(\xi,g):(X,\xi)\in L, g\in \RR\}$
 denote the Jacobi-Dirac structure associated to the Dirac
manifold $(P,L)$. Then $L^c$ naturally has a precontact form as
described in \eqref{precont-form}. The algebroids $\bar{L}$, $L^c$
and $L$ fit into the following diagram (where we denote dimensions
and ranks by superscripts):

\[\xymatrix{  \bar{L}^{n+2} \ar[d]  & (L^c)^{n+1} \ar[r]
\ar[d] & L^n \ar[dl]
\\
Q^{n+1}\ar[r]^{\pi} &P^n&}
\]

The left two Lie algebroids in the diagram are related by the
reduction described in the next proposition:
\begin{prop} \label{algebroid-s1-red}
When $(Q, \bL)$ is a prequantization of Dirac manifold $(P, L)$ we
have $J^{-1}(0)=\bL_0$ (recall that $\bL_0$ was defined
at the end of Section \ref{algoid})
and the isomorphisms of precontact
manifolds and Lie  algebroids,
 \[
\bL//_0 S^1 \cong L^c . \]
\end{prop}
\begin{proof}
 The equality is clear from the characterization of
$J^{-1}(0)$ in eq. \eqref{surj} and from the definition
of $\bL_0$. For the isomorphism notice that $L^c=\pi_{\star} \bar{L}$ (this is
equivalent to saying that $\pi$ is a forward Jacobi-Dirac map) and
apply Prop. \ref{algebroid-red} (which holds because the
assumption $\bar{L}\cap(\g_q,0)\oplus(0,0)=\{0\}$ is satisfied, as
is clear from the definition of $\bL$ in Theorem \ref{thmpreq}).
\end{proof}

In the rest of this subsection we want to see what Lemma \ref{algebroid-s1-red}
says about the objects that integrate the Lie algebroids $\bar{L}$ and $L^c$. We first recall few definitions.

\begin{defi}
A \emph{Lie groupoid} over a manifold $P$
 is  a manifold $\Gamma$ endowed with surjective
  submersions $\bs$,$\bt$ (called source and
 target) to the base manifold $P$, a smooth associative multiplication $m$
 defined on elements $g,h\in \Gamma$ satisfying $\bs(g)=\bt(h)$, an embedding
 of $P$ into $\Gamma$ as the spaces of ``identities'' and a smooth inversion map $\Gamma
 \rightarrow
 \Gamma$ satisfying certain compatibility conditions (see for example   \cite{MW}) \end{defi}
Every Lie algebroid $\Gamma$ has an associated Lie algebroid,
whose  total space is
$ker(\bs_*|_{P})\subset T\Gamma|_P$, with a bracket on sections defined using right invariant
vector fields on $\Gamma$ and   $\bt_*|_{P}$ as anchor. A Lie algebroid $A$ is said to be integrable if there exists a Lie
groupoid whose associated Lie algebroid is isomorphic to $A$; in this case there is a unique (up to isomorphism) source simply connected (s.s.c.)
Lie
groupoid integrating $A$.

 The following two definition are adapted
from  \cite{bcwz} and \cite{IW} respectively  to match up the conventions of \cite{cz} and
\cite{zz}. 
\begin{defi}\label{presymplgroid}
A \emph{presymplectic groupoid} is a Lie groupoid $\Gamma$ over a manifold $P$,
with $\dim \Gamma=2\dim P$,
equipped with a closed 2-form $\Omega_{\Gamma}$   satisfying
$$m^*\Omega_{\Gamma}=pr_1^*\Omega_{\Gamma} +pr_2^*\Omega_{\Gamma}$$
and the non-degeneracy condition
$$\ker \bt_*\cap \ker \bs_* \cap \ker \Omega_{\Gamma}  =\{0\}.$$
\end{defi}
By \cite{bcwz} the Dirac structure on $\Gamma$  given by the graph of $\Omega$
pushes down via $\bs$ to a Dirac structure $L$ on the base $P$ which, as a Lie algebroid, is isomorphic to the Lie algebroid of $\Gamma$. Conversely, if $(P,L)$ is any Dirac manifold, then $L$ (if integrable) integrates to a s.s.c. presymplectic groupoid as above. The latter is unique (up to presymplectic groupoid
automorphism), and will be denoted by $\Gamma_s(P)$ in this paper.\\

Hence presymplectic groupoids are the objects integrating Dirac structures. The objects integrating Jacobi-Dirac structures are the following:
\begin{defi}\label{prect}
A \emph{precontact groupoid} is a Lie groupoid $\Gamma$ over a manifold $Q$,
$\dim \Gamma=2\dim Q+1$,
equipped with a 1-form $\theta_{\Gamma}$ and a function
$f_{\Gamma}$ satisfying
$f_{\Gamma}(gh)=f_{\Gamma}(g)f_{\Gamma}(h)$ and
$$m^*\theta_{\Gamma}=pr_1^*\theta_{\Gamma}pr_2^*f_{\Gamma}+pr_2^*\theta_{\Gamma}$$
and the non-degeneracy condition
$$\ker \bt_*\cap \ker \bs_* \cap \ker \theta_{\Gamma} \cap \ker
d\theta_{\Gamma}=\{0\}.$$
\end{defi}
The 1-form $\theta_{\Gamma}$, viewed as   a Jacobi-Dirac
structure on $\Gamma$,  pushes forward via the source
map to a Jacobi-Dirac structure  on $M$ which
 is isomorphic to the Lie algebroid of $\Gamma$. (The formula for a canonical isomorphism is given in Appendix A). Conversely, if $(Q,\bar{L})$ 
 is any Jacobi-Dirac manifold, then $\bar{L}$ (if integrable) integrates to a s.s.c. unique precontact groupoid as above, which will be denoted by $\Gamma_s(P)$ in this paper. Notice that a Dirac manifold $(P,L)$, in addition to the presymplectic groupoid 
$ \Gamma_s(P)$
 associated as above, also has has an associated precontact
 groupoid $\Gamma_c(P)$ integrating the Jacobi-Dirac structure $L^c$ corresponding 
 to $L$.\\

When the presymplectic groupoid $\Gamma_s(P)$ is prequantizable
its prequantization circle bundle can be view as an ``alternative prequantization space'' for $(P,L)$, because $\Gamma_s(P)$ is the global object that corresponds to the Dirac manifold $(P,L)$. We will see in items (4) and (5) of  Thm. \ref{gen-2} that the
prequantizability and integrability of $(P,L)$ implies that
$\Gamma_s(P)$ is prequantizable, and that the prequantization
bundle $\tilde{\Gamma}_c(P)$ is a groupoid integrating $L^c$, so 
 $A(\tilde{\Gamma}_c(P))\cong L^c$ where ``$A$'' denote the functor that takes the Lie algebroid of a Lie groupoid. (In
the Poisson  case this follows from \cite{cz} and \cite{bcz}).

There is a  canonical  Lie algebroid isomorphism between $\ker
\bs_*|_P \subset T\tilde{\Gamma}_c(P)|_P$ and $L^c$,
  given by 
  Lemma
\ref{idef}. It matches the restriction to $\ker
\bs_*|_P$ of the 1-form on
$\tilde{\Gamma}_c(P)$ and the precontact form  $\theta_{L^C}$ 
on $L^c$ (see eq. \eqref{precont-form})
 at points of $P$   (notice that  at points of
the zero section $P$ the precontact form on $L^c$ is just
$pr^*dt$, i.e. the projection onto the last component).
Similarly the canonical isomorphism between $\ker
\bs_*|_Q$
(where here $\bs$ denotes the source map of $\Gamma_c(Q)$) and $\bar{L}$ 
matches the restriction of the 1-form on $\Gamma_c(Q)$ and $\theta_{\bar{L}}$.
Hence the reduction of Prop. \ref{algebroid-s1-red} matches the 1-forms on the groupoids 
$\Gamma_c(Q)$ and  $\tilde{\Gamma}_c(P)$ at points of the identity sections.

As we will see in the next section, there is an $S^1$ action on
the precontact groupoid
$(\Gamma_c(Q),\theta_{\Gamma},f_{\Gamma})$ of $(Q,\bar{L})$   which is canonically induced by the $S^1$
action on $Q$
 and which hence makes the source map  equivariant
and which respects the 1-form and multiplicative function on the
groupoid. The equivariance makes sure that taking derivatives
along the identity one gets an $S^1$ action on $\ker \bs_*|_Q$ by
vector bundle isomorphism. Further, under the canonical isomorphism (see Lemma
\ref{idef})
 $\ker \bs_*|_Q\cong \bar{L}$,   the $S^1$
action is the natural one described at the beginning of the proof
of Prop. \ref{algebroid-red}, because the $S^1$ action on
$\Gamma_c(Q)$ respects $\bt$,$r_{\Gamma}$ and $\theta_{\Gamma}$.
  We conclude that the $S^1$ action we considered in this
subsection is the infinitesimal version of the $S^1$ action on
$(\Gamma_c(Q),\theta_{\Gamma})$.
We summarize:
\begin{prop}\label{infaction}
The natural $S^1$ action on $Q$ lifts to an action on $A(\Gamma_c (Q))\cong \bar{L}$, whose precontact reduction
is $L^c \cong A(\tilde{\Gamma}_c(P))$, endowed with the Lie algebroid and
precontact structures given by the Lie groupoid $\tilde{\Gamma}_c(P)$.
\end{prop}

In the next section  we will show
that the precontact reduction of $\Gc(Q)$ is isomorphic, both as
precontact manifold and a groupoid, to the s.s.c. precontact groupoid of
$P$, and that $\tilde{\Gamma}_c(Q)$ is a discrete quotient of it.
This means that  precontact reduction commutes
with the Lie algebroid functor:
 
\[ A(\Gamma_c(Q)//_0 S^1)= A(\Gamma_c (Q)) //_0 S^1. \]
Further we also have a correspondence at the intermediate step of
the reduction, namely for the  zero level sets of the moment maps
(see  item $(3)$ of Thm. \ref{redDirac}).

\section{Prequantization and reduction of precontact groupoids}\label{grpic}

In this section we analyze the relation between the groupoids
associated to $(P,L)$ and $(Q,\bL)$, leading to an ``integrated''
version of Proposition \ref{algebroid-s1-red} (i.e. to reduction
of groupoids). In Subsection \ref{poisfd} we will perform the
reduction using finite dimensional arguments, restricting
ourselves for simplicity to the case when $P$ is a Poisson
manifold. If on one hand our finite dimensional proof might appeal
more to geometric intuition, it will not allow to conclude whether
the reduced groupoids we obtain are source simply connected. In
Subsection \ref{dirid},
 for the general case when $P$ is a Dirac manifold,
 we will obtain a complete description of the
reduction using path spaces. We will conclude with two examples.

\subsection{The Poisson case}\label{poisfd}

In this subsection we show our results for Poisson manifold
without using the infinite dimensional path spaces.

We start displaying a simple example, which was  also a motivating
example in \cite{Cr}.
\begin{ep}\label{sympl}
Let $(P,\omega)$ be a simply connected integral symplectic
manifold, and $(Q,\theta)$ a prequantization. We have the
following diagram of groupoids:
\[
\xymatrix{ (Q\times Q\times \RR,-e^{-s}\theta_1+\theta_2,e^{-s})
\ar[d]\ar@<-1ex>[d] & (Q\times_{S^1} Q,[-\theta_1+\theta_2])
\ar[d]\ar@<-1ex>[d] \ar[r] &
(P\times P,-\omega_1+\omega_2) \ar[dl]\ar@<-1ex>[dl] \\
Q \ar[r] & P &  }
\]
The first groupoid is a (usually not s.s.c.) contact groupoid of
$(Q,\theta)$, with coordinate $s$ on the $\RR$ factor. The second
is a contact groupoid of $(P,\omega)$ which is a prequantization
of the third groupoid (the s.s.c. symplectic groupoid of
$(P,\omega)$). The $S^1$ action on $Q$ induces a circle action on
its contact groupoid with moment map given by $\langle J, 1\rangle
=-e^{-s}+1$, so that its zero level set is obtained setting $s=0$,
and dividing by the circle action we obtain exactly the second
groupoid above, i.e. the prequantization of the s.s.c. groupoid of
$(P,\omega)$.
\end{ep}

Let $P$ be a Poisson manifold, consider the Dirac structure $L$ given by the graph of the Poisson bivector,
and assume that $(P,L)$ is
prequantizable and  that it is integrable, in which case it integrates  to a s.s.c symplectic\footnote{This means that 
the 2-form on the presymplectic groupoid intergrating $L$ is non-degenerate.}
 groupoid
$\Gamma_s(P)$.  The prequantizability of $(P,L)$ implies that the
period group of any source fiber of $\Gamma_s(P)$ is contained in
$\ZZ$ (see Section
3.3 of \cite{bcz}, or Theorem \ref{red} below for a straightforward
generalization). This last
condition is equivalent to saying that the symplectic groupoid
$\Gs(P)$ is prequantizable in the sense of \cite{Cr} (see  Prop. 2 in \cite{bcz} or Thm. 3 in \cite{cz}). Its unique
prequantization will be denoted by $\tilde{\Gamma}_c(P)$ and turns
out to be a (usually not s.s.c.) contact\footnote{This means that 
the 1-form on the precontact groupoid satisfies $\theta_{\Gamma}\wedge 
(d\theta_{\Gamma})^{{dim(P)}}\neq 0$ .}
 groupoid of $P$, i.e. it
integrates the Lie algebroid $L^c$. Fix a prequantization $(Q,\bar{L})$ and assume
that the Lie algebroid $\bar{L}$ is integrable; denote by $\Gamma_c(Q)$
the integrating  s.s.c. contact groupoid.
Now,  ``integrating'' the
reduction statements of the last section, we will clarify the
relation between  
$\Gamma_c(Q)$ (the global object attached to the prequantization bundle $Q$) and the prequantization of $\Gamma_s(P)$ (which  
 can be thought of as a different way to prequantize $(P,L)$).
 
The  (smooth)
groupoids we consider  fit into the following diagram; we omitted
$\tilde{\Gamma}_c(P)$, which is just a discrete quotient of the
s.s.c. contact groupoid ${\Gamma}_c(P)$. This diagram corresponds to 
the diagram of  Lie algebroids in Subsection \ref{algpic}, and again we denote
dimensions by superscripts.

\[  
\xymatrix{ \Gamma_c(Q)^{2n+3} \ar[d]\ar@<-1ex>[d] &
\Gamma_c(P)^{2n+1} \ar[d]\ar@<-1ex>[d] \ar[r] &
 \Gamma_s(P)^{2n} \ar[dl]\ar@<-1ex>[dl] \\
Q^{n+1} \ar[r]^{\pi} & P^n &  }
\]
 
\begin{thm}\label{red}
Let $(P,L)$ be an integrable prequantizable Poisson manifold, and
$(Q^{n+1},\bar{L})$ one of its prequantizations as in Subsection \ref{alan}, which we assume
to be integrable. Then:
\begin{itemize}
\item[a)] The s.s.c contact groupoid $\Gamma_c(P)$ of $(P,L)$ is obtained from
the s.s.c. contact groupoid $\Gc(Q)$ of $(Q,\bar{L})$ by $S^1$
contact reduction.
 \item[b)] The prequantization of the s.s.c.
symplectic groupoid $\Gs(P)$ is a discrete quotient of
$\Gamma_c(P)$.
\end{itemize}
\end{thm}

\begin{proof}
$S^1$ acts on $Q$, and it acts also on $TQ\oplus T^*Q$ by the
tangent and cotangent lifts. The $S^1$ action preserves the
subbundle given by the Jacobi-Dirac structure $\bL$, hence we
obtain an $S^1$ action on the Lie algebroid $\bL\rightarrow Q$. The
source simply connected (s.s.c.) contact groupoid $(\Gc(Q),
\theta_{\Gamma},f_{\Gamma})$ of $(Q,\bL)$ is constructed
canonically from the Lie algebroid $\bL$ via the path-space
construction \cite{cf}, so it inherits an $S^1$ action that preserves its
geometric and groupoid structures. In particular the source and
target maps are $S^1$ equivariant, and similarly the
multiplication map $\Gc(Q)_{\bs}\times_{\bt}\Gc(Q) \rightarrow
\Gc(Q)$. Also, the $S^1$ action preserves the contact form, so
there is a moment map $J_{\Gamma}:\Gc(Q)\rightarrow \RR$ by
$J_{\Gamma}(g)=\theta_{\Gamma}(\vG(g))$ where $\vG$ denotes the
infinitesimal generator of the $S^1$ action.
We divide the proof in three steps.\\

\emph{Step 1: $J_{\Gamma}^{-1}(0)$ is a s.s.c. Lie subgroupoid of $\Gc(Q)$.}\\
 We start by showing  that $J_{\Gamma}=1-f_{\Gamma}$; this explicit\footnote{The claim of Step 1 follows even without knowing the
 explicit formula for $J_{\Gamma}$. Indeed
one can show that $J_{\Gamma}^{-1}(0)$ is a subgroupoid by means of the
identity $J_{\Gamma}(gh)=f(h)J_{\Gamma}(g)+J_{\Gamma}(h)$, which is derived
 using the multiplicativity of $\tG$ and the
fact that
 $\vG$ is a multiplicative vector field
(i.e. $\vG(g)\cdot\vG(h)=\vG(gh)$ ; this is just the infinitesimal
version of the statement that the multiplication map is $S^1$
equivariant). Since $J_{\Gamma}^{-1}(0)$ is a smooth wide subgroupoid it
is transverse to the $\bs$ fibers nearby the identity, therefore
its source and target maps are submersions and hence it is
actually a Lie subgroupoid.}
 formula will turn out to be necessary in Step 2.

  To do this we will
use several properties of contact groupoids, for which to refer to
Remark 2.2 in \cite{zz}. The identity $J_{\Gamma}+f_{\Gamma}=1$ is clear
along the identity section $Q$, since $f_{\Gamma}$ is a
multiplicative function and $\vG$ is tangent to $Q$ which is a
Legendrian submanifold of $(\Gc(Q),\tG)$. So to  show that the
statement holds at any point of $\Gc(Q)$ it is enough to show that
$\langle d(f_{\Gamma}+J_{\Gamma}), X_{f_{\Gamma}\bt^*u} \rangle=0$ for
functions $u\in C^{\infty}(Q)$, since hamiltonian vector fields
$X_{f_{\Gamma}\bt^* u}$ span $\ker \bs_*$. The statement follows
by two computations: first
\begin{equation}\label{first}
\begin{split}
\langle
& df_{\Gamma}, X_{f_{\Gamma}\bt^*u} \rangle= \langle df_{\Gamma},
f_{\Gamma}\bt^*uE_{\Gamma}+\Lambda_{\Gamma}d(f_{\Gamma}\bt^*u)\rangle\\
=& f_{\Gamma}\cdot \langle df_{\Gamma},
\Lambda_{\Gamma}d(\bt^*u)\rangle=-f_{\Gamma}\cdot
d(\bt^*u)X_{f_{\Gamma}} =f_{\Gamma}\cdot E(u),
\end{split}
\end{equation}
 where we used twice
$E_{\Gamma}(f_{\Gamma})=0$ and the fact that $\bt$ is a
$-f_{\Gamma}$-Jacobi map. Second,
$$\langle d(\theta_{\Gamma}(v_{\Gamma})), X_{f_{\Gamma}\bt^*u}\rangle =-d\theta_{\Gamma}(v_{\Gamma}, X_{f_{\Gamma}\bt^*u})=
\langle -d (f_{\Gamma}\bt^*u),
(v_{\Gamma}-\theta_{\Gamma}(v_{\Gamma})E_{\Gamma})\rangle
=-f_{\Gamma}\cdot E(u),$$ where we use the fact that
$\cL_{v_{\Gamma}}\theta_{\Gamma}=0$ in the first equality, the
formula $d\theta_{\Gamma}(X_{\phi},w)=-\langle d\phi,w^H \rangle$
valid for any function $\phi$ on a contact groupoid (where $w^H$
is the projection of the tangent vector $w$ to $\ker
\theta_{\Gamma}$ along the Reeb vector field $E_{\Gamma}$) in the
second one, and in the last equality that
$E_{\Gamma}(f_{\Gamma})$,$v_{\Gamma}(f_{\Gamma})$,$\bt_*E_{\Gamma}$
all vanish and that the $S^1$ actions on $\Gamma_c(Q)$ and $Q$ are
intertwined by the target map $\bt$.

Since $f_\Gamma$ is multiplicative, it is clear that
$J_{\Gamma}^{-1}(0)=f_\Gamma^{-1} (1)$ is a subgroupoid.

 Further $J_{\Gamma}^{-1}(0)$ is a smooth submanifold of $\Gc(Q)$:
 by Prop. 3.1.4 in \cite{willett} $g\in \Gc(Q)$ is a singular point of $J_{\Gamma}$ iff
 $\vG(g)$ is a non-zero multiple of $E_{\Gamma}(g)$. Since
 $\tG(E_{\Gamma})=1$   this is never
 the case if $g\in J_{\Gamma}^{-1}(0)$, so $0$ is a regular value of $J_{\Gamma}$.

 To show that $J_{\Gamma}^{-1}(0)$ is a Lie subgroupoid we still need to
 show that its source and target maps are submersions onto $Q$.
 We do so by showing explicitly that $(\ker \bs_*\cap \ker
 d\fG)$ (which along $Q$ will be the Lie algebroid of $J_{\Gamma}^{-1}(0)$) has rank one less than $\ker \bt_*$;
  this is clear since
 by the first equation of Step 1 it is just
$\{X_{f_{\Gamma}\bt^*\pi^*v}: v\in C^{\infty}(P)\}$.

For the proof of the source simply connectedness of the
subgroupoid $J^{-1}_1(0)$  we refer to Thm. \ref{redDirac}.\\

 \emph{Step 2: The contact reduction $J_{\Gamma}^{-1}(0)/S^1$ is the
s.s.c. contact groupoid $\Gamma_c(P)$ of $P$.}

$J_{\Gamma}^{-1}(0)/S^1$ is smooth because the $S^1$ action is  free and
proper, and by contact reduction it
 is a contact manifold,
 so we just have to show that the
Lie groupoid structure descends and is a compatible one.

The $S^1$ equivariance of the source and target maps of $\Gc(Q)$
ensure that source and target descend to maps
$J_{\Gamma}^{-1}(0)/S^1\rightarrow P(=Q/S^1)$. Since
%the $S^1$ equivariance
%of the multiplication on $\Gc(Q)$
%$\vG$ is a multiplicative vector field,
the multiplication on $\Gc(Q)$ is $S^1$ equivariant, the
multiplication on $J_{\Gamma}^{-1}(0)$  induces a multiplication on
$J_{\Gamma}^{-1}(0)/S^1$. It is routine to check this makes
$J_{\Gamma}^{-1}(0)/S^1$ into a groupoid over $P$. Further, since the
source map intertwines the $S^1$ action on $J^{-1}(0)$ and the
free $S^1$ action on the base $Q$, the source fibers of
$J_{\Gamma}^{-1}(0)/S^1$ will be diffeomorphic to the corresponding
source fibers of $J_{\Gamma}^{-1}(0)$, hence we obtain a s.s.c. Lie
groupoid. Since $J_{\Gamma}^{-1}(0)\rightarrow J_{\Gamma}^{-1}(0)/S^1$
 is a surjective submersion,  the $f_\Gamma$-twisted multiplicativity of
 $\theta_\Gamma$ implies that the induced 1-form
$\hat{\theta}_\Gamma$
 is multiplicative, i.e.
 $(J_{\Gamma}^{-1}(0)/S^1,\hat{\theta}_\Gamma ,\hat{f}_{\Gamma})$ is a
 contact groupoid. % \mcomment{Well, in the definition of "contact groupoid" the should be some non-degeneracy too...but we are safe in the $P$=Poisson case}.

%\mcomment{Next argument requires $P$ Poisson}

In order to prove that the above contact groupoid corresponds to
the original Poisson structure $\Lambda_P$ on $P$,
 we have to show that the
source map $\hat{\bs}:J_{\Gamma}^{-1}(0)/S^1\rightarrow P$ is a Jacobi
map (i.e. a forward Jacobi-Dirac map). Consider the diagram
\[
\begin{CD}
J_{\Gamma}^{-1}(0) @>\pi_{J_{\Gamma}}>>  J_{\Gamma}^{-1}(0)/S^1\\
\bs@VVV \hat{\bs}@VVV\\
Q @>\pi>> P.
\end{CD}
\]
We adopt the following short-form notation: for a 1-form $\alpha$,
$L_{\alpha}$ will denote the Jacobi-Dirac structure associated to
$\alpha$ \cite{Wa}. Then for the pullback Jacobi-Dirac structure
we have $i^*L_{\theta_\Gamma}=L_{i^*\theta_\Gamma}$, where $i$ is
the inclusion of $J_{\Gamma}^{-1}(0)$ into $\Gamma_c(Q)$, and the reduced
1-form is recovered as
${\pi_{J_{\Gamma}}}_*i^*L_{\theta_\Gamma}=L_{\hat{\theta}_\Gamma}$. So by
the functoriality of the pushforward, it is enough to show that
$\pi_*\bs_*L_{i^*\theta_\Gamma}$, which by definition is
\begin{eqnarray}\label{push}
\{((\pi \circ\bs)_*Y,f)\oplus (\xi,g): (Y,f)\oplus ((\pi
\circ\bs)^*\xi,g)\in L_{i^*\theta_\Gamma} \},\end{eqnarray}
 equals the Jacobi-Dirac
structure given by $\Lambda_P$. First we determine which tangent
vectors $Y$ to $J_{\Gamma}^{-1}(0)$ and $f\in \RR$ have the property that
$i^*(d\theta_{\Gamma}(Y)+f\theta_{\Gamma})$ annihilates $\ker (\pi
\circ\bs)_*$, which using equation \eqref{first} is equal to
$\{X_{f_{\Gamma}\bt^*\pi^*v}:v\in C^{\infty}(P)\}\oplus \RR
v_{\Gamma}$. A computation similar to those carried out in Step 1
and using the explicit formula $J=1-f_{\Gamma}$ shows that this is
the case when $f=0$ and $\pi_*\bt_*Y=0$, which by a computation
similar to \eqref{first} amounts to $Y\in\{X_{\bs^*\pi^*v}:v\in
C^{\infty}(P)\}\oplus \RR v_{\Gamma}$. These will be exactly the
``$Y$'' and ``$f$'' appearing in  \eqref{push}; a short
computation using the facts that the source map of $\Gamma_c(Q)$
and $\pi$ are Jacobi maps shows that  \eqref{push} equals
$\{(-\Lambda_P \xi,0)\oplus(\xi,g):\xi \in T^*P,g\in \RR\}$, as
was to be shown.
\\

\emph{Step 3: $((J_{\Gamma}^{-1}(0)/S^1)/\ZZ,\hat{\theta}_\Gamma)$ is the
prequantization of the s.s.c. symplectic groupoid ${\Gamma}_s(P)$
of $P$. Here $\ZZ$ acts as a subgroup of
$\RR$ by the flow of the Reeb vector field $\hEG$.}\\
 Consider the action on $J_{\Gamma}^{-1}(0)/S^1$ by its Reeb vector
field $\hat{E}_\Gamma$, which by the contact reduction procedure
is the projection of the Reeb vector field $\EG$ of $\Gc(Q)$ under
$J_{\Gamma}^{-1}(0)\rightarrow J_{\Gamma} ^{-1}(0)/S^1$. 

The $\bt$-image of a $\vG$ orbit is an orbit of the $S^1$ action
on $Q$, since the target map is $S^1$ equivariant. Hence each
$\vG$ orbit meets each $\bt$-fiber at most once. Further each
$\EG$-orbit is contained in a single $\bt$-fiber (since
$\bt_*E_{\Gamma}=0$), so an $\EG$ orbit meets any
 orbit of the $S^1$ action on $\Gc(Q)$ at most once.
Therefore the
 period of an $\EG$ orbit and of the corresponding $\bar{E}_\Gamma$
orbit are equal, and the first period is always an integer number
(because $\bs_*\EG=E_Q$, the generator of the circle action on
$Q$).

Now the we know that the periods of $\bar{E}_{\Gamma}$ are
integers, we can just apply Theorems 2 and 3 of \cite{cz} to prove
our claim.
\end{proof}

\subsection{Path space constructions and the general Dirac case}\label{dirid}

In this subsection we generalize Thm. \ref{red} allowing $P$ to be
a general Dirac manifold, using the explicit description of Lie
groupoids as quotients of path spaces as a powerful tool. The generalization will be presented in Thm. \ref{redDirac} and Thm. \ref{gen-2}.
\\

\begin{defi}Let $\pi:A\rightarrow M$ be a Lie algebroid with anchor $\rho$. 
 The \emph{$A$-path space} $P_a(A)$
consists of all paths $a:[0,1]\rightarrow A$ satisfying 
$\frac{d}{dt}(\pi\circ a)(t)=\rho (a(t))$. 
\end{defi}

There is an equivalence
relation in $P_a A$, called \emph{$A$-homotopy} \cite{cf}.
\begin{defi}
Let $a(t, s)$ be a family of $A$-paths which is $C^2$ in
$s$. Assume that the base paths $\gamma(t,
s):=\pi\circ a(t,s)$ have fixed end points. For a
connection $\nabla$ on $A$, consider the equation
\begin{equation}\label{eq:hom}
\partial_t b-\partial_s a=T_\nabla(a,b),\quad b(0,s)=0.
\end{equation} Here $T_\nabla$ is the torsion of the connection defined
by $T_\nabla(\alpha, \beta)=\nabla_{\rho(\beta)}\alpha-\nabla_{\rho(\alpha)}\beta +
[\alpha,\beta].$ Two paths $a_0=a(0,\cdot)$ and $a_1=a(1,\cdot)$ are
homotopic if the solution $b(t,s)$ satisfies $b(1,s)=0$.
\end{defi}

More geometrically,  for
every Lie algebroid $A$, (notice that tangent bundles are Lie algebroids),  we  associate $A$
a simplicial set $S(A)=[...S_2(A)\Rrightarrow S_1(A) \Rightarrow
S_0(A)]$ with,
\begin{equation} \label{eq:simp-set}
S_i(A) = \hom_{algd}(T\Delta^i , A):= \{ \text{Lie algebroid morphisms}\;T\Delta^i \overset{f}{
\to} A\},
\end{equation}
and face and degeneracy maps $d^n_i: S_n(A) \to S_{n-1}( A) $ and
$s_i^n: S_n(A)\to S_{n+1} (A) $ induced from the natual face and
degeneracy maps $\Delta^n \to \Delta^{n-1} $ and $\Delta^{n}\to
\Delta^{n+1}$.  Here $\Delta^i$ is the $i$-dimensional standard
simplex viewed as a smooth Riemannian manifold with boundary, hence it
is isomorphic to the $i$-dimensional closed ball. Then as explained in
\cite[Section~2]{z:lie2},
\begin{itemize}
\item  it is easy to check that $S_0=M$;
\item  $S_1$ is exactly the $A$-path space $P_a A$;
\item bigons in $S_2$ are exactly the $A$-homotopies in $P_a A$ since a bigon $f: T (d^2_2)^{-1}(Ts^1_0 (T\Delta^0)) \to A$ can
be written as $a(t,s)dt + b(t,s)ds$ over the base map $\gamma(t,s)$
after a suitable choice of parametrization\footnote{We need the one with 
$\gamma(0,s)=x$ and $\gamma(1,s)=y$ for all $s\in[0,1]$. } of
the disk $(d^2_2)^{-1} (s^1_0( \Delta^0))$.
Then we naturally have $b(0,s)=f(0,s)(\frac{\partial}{\partial s})=0$ and
$b(1,s)=f(1,s)(\frac{\partial}{\partial s})=0$.  Moreover the morphism is a
Lie algebroid morphism if and only if $a(t,s)$ and $b(t,s)$ satisfy 
equation  \eqref{eq:hom} which defines the $A$-homotopy. 
\end{itemize}

 The
s.s.c. groupoid of any integrable Lie algebroid $A$ can be constructed
as the quotient of the $A$-path space
 by a foliation $\cF$, whose leaves consists of the $A$-paths that are $A$-homotopic to each other
\cite{cf}. In particular the precontact groupoid $(\Gc(Q), \theta,
f)$ of a Jacobi-Dirac manifold $Q$ can be constructed via the
$A$-path space $P_a (\bL)$, with $\theta$ and $f$ coming from a
corresponding 1-form and function on the path space. We refer to
\cite{cz} \cite{Cr} \cite{IW} and summarize the results in Thm.
\ref{cgpd} below. The advantage of this method is that it can be used
 to generalize
Theorem \ref{red} to the setting of Dirac manifolds (see Theorems \ref{redDirac} and \ref{gen-2}) and that it can be applied to a  general group $G$ action as in \cite{for}.

\begin{thm}\label{cgpd} The s.s.c. precontact groupoid $(\Gc(Q), \theta_\Gamma, f_\Gamma)
$ of an integrable Jacobi-Dirac manifold $(Q, \bL)$ is the
quotient space of the $A$-path space $P_a (\bL)$ by
$A$-homotopies, and $\theta_\Gamma$ and $f_\Gamma$ come from a
1-form $\tilde{\theta}$ and a function $\tilde{f}$ on $P_a(\bL)$.
At the point $a=(a_4, a_3, a_1, a_0)\in P_a (\bL)$, where $(a_4,
a_3, a_1, a_0)$ are components in $TQ \oplus \RR \oplus T^*Q
\oplus \RR$, $\tilde{\theta}$ and $\tilde{f}$ are
\begin{equation} \label{contact-form}
\begin{split}
\tilde{\theta}_{a}(X)= & - \int_0^1 \left\langle e(t) X(t), d
\left(\int_0^1 a_0(t)dt\right) \right\rangle dt + \int_0^1
\left\langle e(t) X(t), pr^*\theta_c \right\rangle dt ,
\\
\tilde{f}(a)=& e(1),\quad \text{with}\; e(t):=e^{\int_0^t - a_3 }
\end{split}
\end{equation}
where $X$ is a tangent vector to $P_a(\bL)$, hence a path itself (parameterized by $t$), and $pr^* \theta_c$ is the pull-back via $pr: \bL\to T^*Q$ of the canonical 1-form on $T^*Q$.  
\end{thm}
\begin{proof}   The equation for $\tilde{f}$ is taken from Prop. 3.5(i) of
\cite{cz}. It is shown there that $\tilde{f}$ descends to the
function $f_\Gamma$ on $\Gc(Q)$.  To get the formula for
$\tilde{\theta}$, we recall from Section 3.4 of \cite{cz} that the
following map $\phi$ is an isomorphism preserving $A$-homotopies:
\[ \phi: P_a (\bL) \times \RR \to P_a(\bL \times_{\psi} \RR), \]
mapping $(a, s)$ with base path $\gamma_1$ to
$\ta:=e^{\gamma_0(t)} a$ with base path $(\gamma_1, \gamma_0)$,
where $\gamma_0:=s-\int^t_0 a_3$. Here $\psi$ is the 1-cocycle on
$\bL$ given by $(X,f)\oplus(\xi,g)\mapsto f$; $\bL \times_{\psi}
\RR$ is the Lie algebroid on $Q\times \RR$ obtained from the Lie algebroid
$\bL$ and the 1-cocycle $\psi$, and it is isomorphic to the Lie
algebroid given by the Dirac structure on $Q\times \RR$ obtained
from the ``Diracization'' of $(Q,\bL)$ (see Section 2.3 in
\cite{IW}).

The correspondence on the level of tangent spaces given by $T\phi$
maps $(\delta \gamma_1, \delta s,\delta a)$ to $(\delta \gamma_1,
\delta \gamma_0, \delta \ta)$ and satisfies
\[
\begin{split}
\delta \gamma_0 &=\delta s - \int_0^t  a_3  , \\
\delta \ta_1   &= e^{\gamma_0}( \delta a_1 + (\delta s- \int_0^t
 \delta a_3   ) a_1), \\
\delta \ta_0 &= e^{\gamma_0} (\delta a_0+ (\delta s-\int_0^t
\delta a_3  ) a_0 ).
\end{split}
\]
We identify $\bL\times_{\psi} \RR$ with the Dirac structure on
$Q\times \RR$ given by the Diracization of $(Q,\bar{L})$. Then on the whole space $P
(\bL\times_{\psi} \RR)$ of paths in $\bL\times_{\psi} \RR$ there
is a symplectic form $\omega$ coming from integrating the
pull-back of the canonical symplectic form on $T^*(Q\times \RR)$
(see Section 5 in \cite{bcwz}). This form restricted to the
$A$-path space $P_a (\bL \times_{\psi} \RR )$ is homogeneous
w.r.t. the $\RR$ component, i.e. $\varphi_s \omega = e^s \omega$,
where $\varphi_s$ is the flow of $\frac{\partial}{\partial s}$
with $s$ the coordinate of $\RR$. This is because $\varphi_s$ acts
on vector fields $\delta \ta_1$ and $\delta \ta_0$ by rescaling by
an $e^s$ factor as the formula of $T\phi$ and $\gamma_0$ show.
This homogeneity survives the quotient to groupoids as shown in
\cite{cz}. Therefore $\theta_\Gamma$ comes from the 1-form
$\tilde{\theta}$ whose associated homogeneous symplectic form is
$\omega$, i.e. $\tilde{\theta} = -i^*_0
i(\frac{\partial}{\partial s}) \omega$. With a straightforward
calculation and the formula of $T\phi$, we have the formula for
$\tilde{\theta}$ in \eqref{contact-form}.
\end{proof}
\begin{remark}

 The formula
for $\tilde{\theta}$ is a generalization of Theorem 4.2 in
\cite{Cr} in the case $\bL$ that comes from a Dirac structure. To
get the formula of the 1-form there up to sign\footnote{In
\cite{Cr} 1-forms on contact groupoids are so that the target map
is a Jacobi map, whereas here we adopt the convention (as in
\cite{zz}) that the source map be Jacobi.}, one just has to put
$e(t)=1$ which corresponds to the case that $a_3=0$.
\end{remark}

In Lemma \ref{pullback} we constructed a Lie algebroid structure
on $\pi^* A$, the pull back via $\pi: Q\to P$ of any Lie algebroid
$A$ on $P$, provided there is a flat $A$-connection $\tilde{D}$ on
the vector bundle $K$ corresponding to the principal bundle $Q$. ($\pi^*A$ turned out to
be the transformation algebroid w.r.t. the action by the flat
connection). Now we show some functorial property of algebroid paths in  $\pi^*A$. Later in this section we will apply them to $A=L^c$,
for $\pi^*L^c$ is identified   with a Lie subalgebroid of $\bar{L}$ (Thm. \ref{summ}),
whose integrating groupoid we can describe in term of $A$-paths (Thm. \ref{cgpd}).

\begin{lemma}\label{homotopy-func}
An $A$-path $a$ in $A$  can be lifted to an $A$-path in $\pi^*A$.
The same is true for $A$-homotopies. In other words,  in the
following diagram (for $n=1,2$),
\[
\xymatrix{ T \Delta^n \ar[rrd] \ar[d] \ar[dr]^{ f} & &
\\
\Delta^n \ar[dr]^{ f_0} \ar[drr] & \pi^* A \ar[r] \ar[d] & A \ar[d] \\
& Q \ar[r]^{\pi} & P }
\]
any Lie algebroid morphism  $f: T\Delta^n \to A$ lifts to a Lie
algebroid morphism from $T\Delta^n$ to $\pi^* A$.
\end{lemma}
\begin{proof}
 Let $\gamma$ be the base path of an $A$-path $a$, and let
 $\tilde{\gamma}$ be the parallel translation along $a$ of some
 $\tg(0)\in \pi^{-1}(
\gamma(0))$ as in the proof of Lemma \ref{pullback}.  Denoting by
$\pi^*a$ the lift of $a$ to $\pi^*A$ with base path
$\tilde{\gamma}$, we have $\rho(\pi^* a) =
h_Q(a(\gamma(t)),\tilde{\gamma}(t))=d/dt (\tg)$, with $\rho$ the
anchor of $\pi^*A$ (see equation \eqref{lifth}). That is, $\pi^*
a$ is an $A$-path in $\pi^* A$ over $\tg$. The lifting of $a$ is
not unique. In fact it is determined by the choice of a point in
$\pi^{-1}( \gamma(0))$ as initial value.

Now we prove the same statement for $A$-homotopies. Suppose
$a(\epsilon, t)$ is an $A$-homotopy over $\gamma(\epsilon, t)$,
i.e. there exist $A$-paths (w.r.t. parameter $\epsilon$)
$b(\epsilon, t)$ also over $\gamma$ satisfying
\begin{equation} \label{a-homotopy}
\partial_t b
 -\partial_{\epsilon}  a
=\nabla_{\rho (b)} a- \nabla_{\rho (a)} b + [a, b],
\end{equation}
and the boundary condition $b(\epsilon, 0)=b(\epsilon, 1)=0$, for
any choice of connection $\nabla$ on $TP$.
 As
above, we can lift $\gamma$ to $\tg(\epsilon, t)$. In fact, once
we choose $\tg(0,0)$, we can use $\tg(0,0)$ to obtain the lift
$\tg(\epsilon, 0)$ and then $\tg(\epsilon, t)$. (The lift does not
depend on whether we lift  $\gamma(\epsilon,0)$ or $\gamma(0,t)$
first,
  because the connection
 $\tilde{D}$ is flat). Then $\pi^* a$ and $\pi^* b$ are $A$-paths
over $\tg$ w.r.t. parameters $t$ and $\epsilon$ respectively.
Moreover, we choose a connection $\tn$ on $Q$ induced from the
connection $\nabla$ on $P$ such that $\tn_{X^H} Y^H = (\nabla_X
Y)^H$, $\tn_{X^H} E = 0$, $\tn_E Y^H =0$ and $\tn_E E= 0$, where
the superscript $H$ denotes the horizontal lift with respect to
some connection we fix on the circle bundle $\pi:Q \rightarrow P$.
(Since $E(\pi^* f) =0$ and $X^H(\pi^*f)=X(f)$ these requirements
are consistent. In fact, the connection $\tn$ on $TQ=\pi^*TP\oplus
\RR E$ is just the sum of the pullback connection on $\pi^*TP$ and
of the trivial connection). Now we will prove that $\pi^* a$ and
$\pi^* b$ satisfy \eqref{a-homotopy} w.r.t. $\tn$. Notice that
$\langle \pi^* \eta, \tn_E X \rangle =0$ for all vector fields
$X$, so we have
\[ \tn_E \pi^* \eta =0, \quad \tn_{(\frac{\partial}{\partial \epsilon }
\gamma)^H}\pi^* \eta = \pi^*(\nabla_{\frac{\partial}{\partial
\epsilon } \gamma} \eta). \] Therefore
$\tn_{\frac{\partial}{\partial \epsilon } \tg} \pi^* \eta = \pi^*
(\nabla_{\frac{\partial}{\partial \epsilon } \gamma} \eta)$. So
$\partial_\epsilon \pi^* a = \pi^*(\partial_\epsilon a)$. The same
is true for $\pi^* b$. Moreover, since
 $\rho (\pi^* a) = (\rho
(a))^H +  \langle \tilde{\beta},a \rangle E$ (upon writing
$\tilde{D}$ as in equation \eqref{deco}
 and denoting by $^H$ the horizontal lift w.r.t. $\ker
\sigma$), similarly we have $\tn_{\rho (\pi^* a)} \pi^* b =
\pi^*(\nabla_{\rho(a)} b)$ as well as the analog term obtained
switching $a$ and $b$. By the definition of Lie bracket on
$\pi^*A$, we also have $[\pi^* a, \pi^* b]=\pi^*([a, b])$.
Therefore $a$, $b$ satisfying \eqref{a-homotopy} implies that the
same equation holds for $\pi^*a $ and $\pi^*b$. The boundary
condition $\pi^*b(\epsilon, 0)=\pi^*b(\epsilon, 1)=0$ is obvious.
Hence, $\pi^* a$ is an $A$-homotopy in $\pi^* A$.
\end{proof}
\begin{remark}\label{rk:full}
We claim that all the $A$-paths and $A$-homotopies in $\pi^*A$ are
of the form $\pi^* a$. Indeed consider a $\pi^*A$ path $\hat{a}$
over a base path $\hat{\gamma}$, i.e. $\rho
(\hat{a}(t))=\frac{d}{dt}\hat{\gamma}(t)$.
 Let $\gamma:=\pi \circ \hat{\gamma}$ and let $a(t)$ be equal to
$\hat{a}(t)$, seen as an element of $A_{\gamma(t)}$. The
commutativity of
\[
\begin{CD}
\pi^*A @>h_Q=\rho >>  TQ\\
@VVV \pi_*@VVV\\
A @>\rho_A >> TP
\end{CD}
\]
implies that $a$ is an $A$-path over $\gamma$. Further, the
horizontal lift of $a$ starting at $\hat{\gamma}(0)$ satisfies by
definition $\frac{d}{dt}\tilde{\gamma}(t)=
h_Q(a(\gamma(t)),\tilde{\gamma}(t))$, so it coincides with
$\hat{\gamma}$. The same holds  for $A$-homotopies.
\end{remark}

The next theorem generalizes  Thm. \ref{red}a).

\begin{thm}\label{redDirac}Let $(P, L)$ be an integrable prequantizable Dirac manifold and $(Q, \bL)$ one of its
prequantization. We use the notation $[\cdot]_A$ to denote $A$-homotopy classes
in the Lie algebroid $A$. Then we have the following results:
\begin{enumerate}
\item there is an $S^1$ action on the precontact groupoid $\Gamma_c(Q)$ with moment map $J_{\Gamma} = 1- f_\Gamma$;
\item $J_{\Gamma}^{-1}(0)$ is a source connected and simply connected subgroupoid  of $\Gamma_c(Q)$ and is isomorphic to the action groupoid $\Gamma_c (P) \ltimes Q \rightrightarrows
Q$.
\item In terms of path spaces,
\[J_{\Gamma}^{-1}(0)=\{ [\pi^* a]_{\bL} \}= \{ [\pi^* a]_{\bL_0} \},\]
where $a$ is an $A$-path in $ L^c$ and $\pi^* a$ is defined as in
Lemma \ref{homotopy-func} (we identify $\pi^*L^c$ with
$\bL_0\subset \bL$ as in Thm. \ref{summ}). Hence   the
Lie algebroid of $J_{\Gamma}^{-1}(0)$ is $\bar{L}_0=J^{-1}(0)$ (see Prop.
\ref{algebroid-s1-red}).
\item the precontact reduction $\Gamma_c(Q)//_0 S^1$ is isomorphic to the s.s.c. contact groupoid $\Gc(P)$ via
the inverse of the following map
\[ p: [a]_{L^c} \mapsto [\pi^* a]_{\bL, S^1},
\]where  $[\cdot]_{\bL, S^1}$ denotes $S^1$ equivalence classes of $[\cdot]_{\bL}$.
\end{enumerate}
\end{thm}
 
\begin{remark} The isomorphism $p$ gives the same contact groupoid
structure on $\Gamma_c(Q)//_0 S^1$ as in Theorem \ref{red} in the
case when $P$ is Poisson.
\end{remark}
\begin{proof}
$1)$ The definition of the $S^1$ action is the same as in Theorem
\ref{red}.  $J_{\Gamma}$ is defined by
$J_{\Gamma}(g)=\theta_\Gamma(v_\Gamma(g))$, where $v_\Gamma$ is induced
by the $S^1$ action on $Q$ hence on $\bL$. More explicitly, $T(P_a
(\bL))$ is a subspace of the space of  paths in $T\bL$. If we take
a connection $\nabla$ on $Q$, then $T\bL$  decomposes as $TQ
\oplus \bL$. At $(a_4, a_3, a_1, a_0) \in P_a(\bL)$ the
infinitesimal $S^1$ action $\tv$ on the path space is $
\tv=(E(\gamma(t)),*,*,*, 0)$.
%term in $T^*Q$ such that $\langle *, v\rangle =\langle -\nabla_v E, \xi\rangle $ when the connection is trivial.
So
\[ J_{\Gamma}([a])= \tilde{\theta}_{a} (\tv) = \int_0^1(\langle a_1(t),
E\rangle e^{-\int^t_0 \langle a_1, E \rangle dt})
dt=- \int_0^1 d (e^{-\int_0^t \langle a_1 , E\rangle dt})=1-  f_\Gamma.\]\\

$2)$ By $1)$ $J_{\Gamma}^{-1}(0)= f_\Gamma^{-1} (1)$. Since $f_\Gamma$ is
multiplicative, it is clear that $f_\Gamma^{-1} (1)$ is a
subgroupoid. Moreover  using Thm. \ref{cgpd} we see that
$f_\Gamma^{-1}(1)$ is made up by paths $a=(a_4, a_3, a_1, a_0)$
such that
\begin{equation}\label{j1} \int_0^1\langle a_1(t), E\rangle dt=0.
\end{equation} Notice that this are not exactly the same as $A$-paths in
$\bL_0$, which are the $A$-paths such that $\langle a_1(t) , E
\rangle \equiv 0$ for all $t\in [0, 1]$ (see Thm. \ref{summ}).

Now we show that $J_{\Gamma}^{-1}(0)$  is source connected. Take $g\in
\bs^{-1}(x)$, and choose an $A$-path $a(t)$ representing $g$ over
a base path $\gamma(t):I\rightarrow Q$. We will connect $g$ to $x$
within $J_{\Gamma}^{-1}(0)\cap \bs^{-1}(x)$ in two steps: first we deform
$g$ to some other point $h$ which can be represented by an
$A$-path in $\bL_0$; then we ``linearly shrink'' $h$ to $x$.

Suppose the vector bundle $\bL$ is trivial on  a neighborhood $U$
of the image of $\gamma$ in $Q$. Choose a frame
$Y_{0},\dots,Y_{\dim Q}$ for $\bL|_U$, with the property that
$Y_0=(-A^H,1)\oplus(\sigma-\pi^*\alpha,0)$ (with $\sigma$, $A$ and
$\alpha$ as in Thm. \ref{thmpreq}) and that all other $Y_i$
satisfy $\langle a_1,E \rangle=0$. In this frame,
$a(t)=\sum_{i=0}^{\dim Q}p_i(t) Y_i|_{\gamma(t)}$ for some
time-dependent coefficients $p_i(t)$. Define the following section
of $\bL|_U$:
$Y_{t,\epsilon}=(1-\epsilon)p_0(t)Y_0+\sum_{i=1}^{\dim Q}p_i(t)
Y_i$.  Define a deformation  $\gamma(\epsilon, t)$ of $\gamma(t)$
by
\[ \frac{d}{dt} \gamma(\epsilon, t) =\rho(Y_{t, \epsilon}), \quad
\gamma(\epsilon, 0)=x,\] where $\rho$ is the anchor of $\bL$ (one
might have to extend $U$ to make $\gamma(\epsilon, t) \in U$ for
$t\in [0, 1]$). Let $a(\epsilon,
t):=Y_{t,\epsilon}|_{\gamma(\epsilon,t)}$. For each $\epsilon$  it
is an $A$-path by construction, and $a(0, t)=a(t)$. Using  $g\in
J_{\Gamma}^{-1}(0)$ (so that $\int_I p_0(t) dt=0$) we have
$$\int_0^1\langle a_1( \epsilon, t),E \rangle
dt= \int_0^1\langle (1-\epsilon)p_0(t)Y_0+\sum_{i=1}^{\dim
Q}p_i(t) Y_i, (E, 0, 0, 0) \rangle_- dt =(1-\epsilon) \int_I
p_0(t) dt=0,$$ so $[a(\epsilon, \cdot)]$ lies in $J_{\Gamma}^{-1}(0)$.
Notice that $a(1,t)$ satisfies $\langle a_1(1, t),E \rangle \equiv
0$ for all $t$; hence an $A$-path in $\bL_0$. We denote $h:=[a(1,
t)]$ and define a continuous map $pr: P_a(\bL|_U) \to P_a(\bL_0
|_U)$ by $a(t)\mapsto a(1, t)$.

Then we can shrink linearly $a(1,t)$ to the zero path, via
$a^{\delta}(1, t):=\delta a(1, \delta t)$ which is an $A$-path
over $\gamma(1, \delta t)$. Taking equivalence classes we obtain a
path from $h$ to $x$, which moreover lies in $J_{\Gamma}^{-1}(0)$ because
$\langle a_1(1, t) ,E \rangle \equiv 0$.
 
Now we show that $J_{\Gamma}^{-1}(0)$ is source simply connected. If
there is a loop $g(s)=[a(1, s, t)]$   in a source fibre of
$J_{\Gamma}^{-1}(0)$, then $g(s)$ can shrink to $x:=\bs(g(s))$ inside the
big (s.s.c.!) groupoid $\Gc(Q)$ via $g(\epsilon, s)=[a(\epsilon,
s, t)]$. We can assume $a(\epsilon, s, t)=sa(\epsilon,1, st)$.
This is easy to realize since we can simply take $a(\epsilon, s,
t)= g(\epsilon, st)^{-1} d/dt(g(\epsilon, st))$. Then the $a(i, 1,
\cdot)$'s are $A$-paths in $\bL_0$ for $i=0, 1$. This is because
both $g(s)$ and $x$ are paths in $J_{\Gamma}^{-1}(0)$ which implies
$\int_0^1 sa(i, 1, st)=0$ for all $s\in[0,1]$. Moreover the base
paths $\gamma(\epsilon, s, t)$ form an embedded disk (one can
assume that the deformation $g(\epsilon, s)$ has no
self-intersections) in $Q$. So we can take a simply connected open
set (for example a tubular  neighborhood of this disk) $U\subset
Q$ containing $\gamma(\epsilon, s, t)$. Then $L|_U$ is trivial.
Therefore there is a continuous map $pr$ such that
$\bar{a}(\epsilon, 1, \cdot) =pr (a(\epsilon, 1, \cdot))$ is an
$A$-path in $\bL_0$ and $\bar{a}(1, 1, \cdot) = a(1,1, \cdot)$.
Then we can shrink $g(s)=\bar{g}(1, s)$ to $x=\bar{g}(0, s)$ via
 \[\bar{g}(\epsilon, s):=[s\bar{a}(\epsilon,1,  st)], \]
which is inside of $J_{\Gamma}^{-1}(0)$ since $\langle
\bar{a}_1(\epsilon, 1, t), E \rangle \equiv 0$.
\\

 $3)$ To show that
$J_{\Gamma}^{-1}(0)= \{ [\pi^* a]_{\bL} \}$, we just have to show that an
$A$-path in $\bL$ satisfying \eqref{j1} is $A$-homotopic
(equivalent) to an $A$-path lying contained in  $\bL_0$. Since
$J_{\Gamma}^{-1}(0)$ has connected source fibres, given a point $g=[a]$
in $J_{\Gamma}^{-1}(0)$, there is a path $g(t)$ connecting $g$ to
$\bs(g)$ lying in $J_{\Gamma}^{-1}(0)$. Differentiating $g(t)$ we get an
$A$-path $b(t)= g(t)^{-1} \dot{g}(t)$ which is $A$-homotopic to
$a$ and $s b(s t)$ represents the point $g(s t) \in J^{-1}(0)$.
Therefore $\int_0^1 \langle s b_1 (s t), E \rangle dt=0$, for all
$s \in [0, 1]$. Hence $\langle b_1(t), E\rangle \equiv 0$ for all
$t\in [0, 1]$, i.e. $b$ is a path in $\bL_0$.

To further show that $J_{\Gamma}^{-1}(0)=\{ [\pi^* a]_{\bL_0} \}$, we
only have to show that if two $A$-paths in $\bL_0$ are
$A$-homotopic in $\bL$ then they are also $A$-homotopic in
$\bL_0$. Let $a(1, \cdot)$ and $a(0, \cdot)$ be two $A$-paths in
$\bL_0$, $A$-homotopic in $\bL$ and representing an element $g\in
J_{\Gamma}^{-1}(0)$. Integrate $sa(i, st)$ to get $g(i, t)$ for $i=0, 1$.
Namely we have $sa(i, st)=g(i, s)^{-1} \frac{d}{dt}|_{t=s} g(i,
t)$.
%Integrating $a(\epsilon, t)$, we get a groupoid homotopy $g(\epsilon, s)$ lying in $\bs^{-1}(x)$ where $x=\bs([a(0, t)]$ such that $g^{-1}(\epsilon, st) d/dt (g(\epsilon, st)) = s a(\epsilon, st)$, namely $g(\epsilon, s)=[sa(\epsilon, st)]$ (especially $g(1,1)=[a(1, t)]=[a(0, t)]$). $g(\epsilon, s)$ might not all lie in $J_{\Gamma}^{-1}(0)$.
Then $g(i, t)$ are two paths connecting $g$ and $x:=\bs(g)$ lying
in the subgroupoid $J_{\Gamma}^{-1}(0)$ since $a(i, t)$ are paths in
$\bL_0$. Since the source fibre of $J_{\Gamma}^{-1}(0)$ is simply
connected, there is a homotopy $g(\epsilon, t) \in J_{\Gamma}^{-1}(0)$
linking $g(0, t)$ and $g(1, t)$. So $sa(\epsilon, st):=
g(\epsilon, s)^{-1} \frac{d}{dt}|_{t=s} g(\epsilon, t)$ is an
$A$-path in the variable $t$ representing the element $g(\epsilon,
s)\in J_{\Gamma}^{-1}(0)$ for every fixed $s$. Hence $sa(\epsilon, st)$
satisfies \eqref{j1} for every $s\in [0, 1]$. Therefore $\langle
a_1(\epsilon, t), E \rangle \equiv 0$. Then $a(\epsilon, t)
\subset \bL_0$ is an $A$-homotopy between $a(0, t)$ and $a(1, t)$.

Therefore $J_{\Gamma}^{-1}(0)$ is the s.s.c. Lie groupoid integrating
$J^{-1}(0)=\bL_0$.

$4)$ First of all,  given an $A$-path $a$ of $L^c$ over the base
path $\gamma$ and a point $\tg(0)$ over $\gamma(0)$ in $Q$, we
lift it to an $A$-path $\pi^*a$ of $\bL$ as described in Lemma
\ref{homotopy-func}. By the same lemma, we see that ($L^c$)
$A$-homotopic $A$-paths in $L^c$ lift to ($\bL_0$) $A$-homotopic
$A$-paths in $ \pi^* L^c \cong \bL_0\subset \bL$, so the map $p$
is well defined  Different choices of $\tg(0)$ give exactly the
$S^1$ orbit of (some choice of) $[\pi^*a]_{\bL}$. Surjectivity of
the map $p$ follows from the statement about $A$-paths in Remark
\ref{rk:full}. Injectivity follows from the fact that $\{ [\pi^*
a]_{\bL} \}=\{[\pi^* a]_{\bL_0}\}$ in $3)$ and the statement about
$A$-homotopies in Remark \ref{rk:full}.
\end{proof}

We saw in Subsection \ref{algpic}
that, given any integrable Dirac manifold $(P, L)$, there are two
groupoids attached to it. One is the presymplectic groupoid
$\Gamma_s(P)$ integrating $L$; the other is the precontact
groupoid $\Gamma_c(P)$ integrating $L^c$. In the non-integrable
case, these two groupoids still exist as stacky groupoids carrying
the same geometric structures (presymplectic and precontact)
\cite{tz}. In this paper, to simplify the treatment, we view them
as topological groupoids carrying the same name and when the
topological groupoids are smooth manifolds they have additional
presymplectic and precontact structures. Item (4) of the following theorem generalizes
 Thm. \ref{red}b). The other items
generalize from the Poisson case to the Dirac case Theorem 2 and 3 in \cite{cz} and a result in \cite{bcz}.
\begin{thm}\label{gen-2}
For a Dirac manifold $(P, L)$, there is a short exact sequence of
topological groupoids
\[1\to \cG \to \Gamma_c(P) \overset{\tau}{\to} \Gamma_s(P) \to 1,\]
where $\cG$ is the quotient of the trivial groupoid $\RR\times P$
by a group bundle $ \cP$ over $P$ defined by
\[
\begin{split}
\cP_x:= &\{ \int_{[\gamma]} \omega_{F}: [\gamma]\in \pi_2(F,
x)\; \text{and}\; \gamma\; \text{is the base of an} \\
& \text{$A$-homotopy between paths representing $1_x$ in $L$.}\},
\end{split}
\] with $F$ the
presymplectic leaf passing through $x\in P$ and $\omega_{F}$ the
presymplectic form on $F$. In the case that $(P, L)$ is integrable
as a Dirac manifold, then
\begin{enumerate}
\item the presymplectic form $\Omega$ on $\Gamma_s(P)$ is related to
the precontact form $\theta$ on $\Gamma_c(P)$ by
\[ \tau^* d \theta = \Omega, \]
and the infinitesimal action $R$ of $\RR$ on $\Gamma_c(P)$ via
$\RR \times P \to \cG$  satisfies
\[ \cL_R \theta =0, \quad i(R)\theta =1. \]
\item  $R$ is the left invariant vector field
extending the section $(0,0)\oplus(0, -1)$ of $L^c\subset
\cE^1(P)$ as in Cor. \ref{computation};
\item the group $\cP_x$ is generated by the periods of $R$;
\item $\Gamma_s(P)$ is prequantizable iff $\cP \subset P\times \ZZ$;
in this case its prequantization is $\Gamma_c(P)/\ZZ$, where  
$\ZZ$ acts on $\Gamma_c(P)$  as a subgroup of $\RR$.
\item If $P$ is prequantizable as a Dirac manifold, then $\Gamma_s(P)$ is
prequantizable.
\end{enumerate}
\end{thm}
\begin{proof}
The proof of $(1)$ and $(4)$ is the same as Section 4 of
\cite{cz}. One only has to replace  the Poisson bivector $\pi$ by
$\Upsilon$ and the leaf-wise symplectic form of $\pi$ by
$\omega_F$. $(3)$ is clear since $R$ generates the $\RR$ action
and $\cG=\RR/\cP$.

For $(2)$, we identify $(0,0)\oplus(0,-1)$ with a section of $\ker
\bt_*$
 using Lemma
\ref{idef} and then extend it to a left invariant vector field on
$J^{-1}(0)/S^1$. Using Cor. \ref{computation} we see that the
resulting vector field is killed by $\bs_*$, $\bt_*$ and
$d\theta_{\Gamma}$ and that it pairs to 1 with $\theta_{\Gamma}$,
so by the ``non-degeneracy'' condition in Def. \ref{prect}
 it must be equal to $R$.

For $(5)$, if $P$ is prequantizable as a Dirac manifold, then
$\Upsilon = \rho^* \Omega+d_L \beta$ for some integral form
$\Omega$ on $P$ and $\beta\in \Gamma(L^*)$. Suppose $f= a
d\epsilon + b dt$ is a Lie algebroid homomorphism from the tangent
bundle $T\square$ of a square $[0,1]\times [0,1]$ to $L$ over the
base map $\gamma: \square \to P$, i.e. $a(\epsilon, t)$ is an
$A$-homotopy over $\gamma$ via $b(\epsilon, t)$ as in
\eqref{a-homotopy}. Denoting by $\omega_F$ the presymplectic form
of the leaf $F$ in which $\gamma(\square)$ lies, we have (see also
Sect. 3.3 of \cite{bcz}),
\[
\begin{split}
\int_\gamma \omega_F &= \int_{\square} \omega_F (\frac{\partial \gamma}{\partial t}, \frac{\partial\gamma}{\partial \epsilon}) = \int_{\square} \langle a d\epsilon, b d t \rangle_- = \int_{\square} f^* \Upsilon \\
&=\int_{\square} f^* (\rho^* \Omega+d_L\beta)= \int_{\square} f^*
(\rho^* \Omega)= \int_{\square} \gamma^* \omega =\int_{\gamma}
\omega \in\ZZ
\end{split}
\]
where we used $\Upsilon=\rho^*\omega_F$ in the second equation and
$f^*d_L\beta=d_{dR}(f^*\beta)$ in the fifth.
\end{proof}

\subsection{Two examples}\label{anex}

We present two explicit examples for Thm. \ref{red},
\ref{redDirac} and \ref{gen-2}.

The first one generalizes Example \ref{sympl}.
\begin{ep}\label{sympl2}
Let $(P,\omega)$ be an integral symplectic manifold (non
necessarily simply connected), and $(Q,\theta)$ a prequantization.
The s.s.c. contact groupoid of $(Q,\theta)$ is
$(\bar{Q}\times_{\pi_1(Q)} \bar{Q}\times
\RR,-e^{-s}\theta_1+\theta_2,e^{-s})$ where $\bar{Q}$ denotes the
universal cover of $Q$. As in Example \ref{sympl} the moment map
is given by $J_{\Gamma}=-e^{-s}+1$ and the reduced manifold at zero is
$((\bar{Q}\times_{\pi_1(Q)}\bar{Q})/{S^1}, [-\theta_1+\theta_2])$,
where $\pi_1(Q)$ acts diagonally and the diagonal $S^1$ action is
realized by following the Reeb vector field on $\bar{Q}$.

Notice that the Reeb vector field of
$(\bar{Q}\times_{\pi_1(Q)}\bar{Q})/{S^1}$ is the Reeb vector field
of the second copy of $\bar{Q}$. Dividing $\bar{Q}$ by
$\ZZ\subset\text{(Flow of Reeb v.f.)}$ is the same as dividing by
the $\pi_1(\tilde{Q})$ action on $\bar{Q}$, where $\tQ$ is the
 pullback of $Q\rightarrow P$ via the universal covering $\tP\rightarrow
 P$. To see this use that $\pi_1(\tilde{Q})$
is generated by any of its Reeb orbits (look at the long exact
sequence corresponding to $S^1\rightarrow \tilde{Q} \rightarrow
\tilde{P}$), and that the Reeb vector field of $\bar{Q}$ is
obtained lifting the one on $\tilde{Q}$. Also notice that
$\pi_1(\tilde{Q})$ embeds into $\pi_1({Q})$ (as the subgroup
generated by the Reeb orbits of $Q$) and that the quotient by the
embedded image is isomorphic to $\pi_1(P)$,
 by the long exact sequence for $S^1\rightarrow Q \rightarrow P$.
 So the quotient of $(\bar{Q}\times_{\pi_1(Q)}\bar{Q})/{S^1}$
 by the $\pi_1(\tilde{Q})$ action on the second factor
 is
 $(\tilde{Q}\times_{\pi_1({P})}\tilde{Q})/S^1$
where we used $\bar{Q}/\pi_1(\tilde{Q})=\tilde{Q}$ on each factor.
This groupoid, together with the induced 1-form
$[-\theta_1+\theta_2]$,
 is clearly the prequantization of the s.s.c. symplectic groupoid
$(\tilde{P}\times_{\pi_1(P)}\tilde{P},-\omega_1+\omega_2)$ of
$(P,\omega)$.
\end{ep}

In the second example we consider a Lie algebra $\g$. Its dual
$\g^*$ is endowed with a linear Poisson structure $\Lambda$,
called Lie-Poisson structure, and the Euler vector field $A$
satisfies $\Lambda=-d_{\Lambda}A$ where $d_{\Lambda}$ is the Poisson cohomology differential. So the prequantization
condition \eqref{cond0} for $(\g^*,\Lambda)$ is satisfied, with
$\Omega=0$ and $\beta=A$. We display the contact groupoid
integrating the induced prequantization $(Q,\bL)$ for the simple
case that $\g$ be one dimensional; then we show that (a discrete
quotient of) the $S^1$ contact reduction of this groupoid is the
prequantization of the symplectic groupoid of $\g^*$.

\begin{ep}\label{1dim}
Let $\g=\RR$ be the one-dimensional Lie algebra. We claim that the
prequantization $Q=S^1\times \g^*$ of $\g^*$ as above has as a
s.s.c. contact groupoid ${\Gamma}_c(Q)$ the quotient of
 \begin{equation}\label{groidgstar}  (\RR^5,xd\epsilon-e^td\theta_1+d\theta_2,e^t)\end{equation}
by the diagonal $\ZZ$ action on the variables
$(\theta_1,\theta_2)$.
%Here $\tilde{Q}\cong \RR^2$ is the
%universal cover of $Q$, and
Here the coordinates on the five factors of $\RR^5$ are
$(\theta_1,t,\epsilon,\theta_2,x)$. The groupoid structure is the
product of the following three groupoids: $\RR\times
\RR=\{(\theta_1,\theta_2)\}$ the pair groupoid; $\RR\times
\RR=\{(t,x)\}$ the action groupoid given by the flow of the vector
field $-x
\partial_x$ on $\RR$, i.e. $(t',e^{-t}x)\cdot(t,x)=(t'+t,x)$; and $\RR=\{\epsilon\}$
the group.

To see this, first determine the prequantization of
$(\g^*,\Lambda)$: it is $Q=S^1\times \RR$ with Jacobi structure
$(E\wedge x
\partial_x,E)$, where $E=\partial_{\theta}$ is the infinitesimal generator of the
circle action and $x
\partial_x$ is just the Euler vector field on $\g^*$ (see \cite{cmdl}). This Jacobi
manifold has two open leaves, and we first focus on one of them,
say $Q_+=S^1\times \RR_+$. This is a locally conformal symplectic
leaf, with structure $(d\theta\wedge \frac{dx}{x},\frac{dx}{x})$.

We determine the s.s.c contact groupoid $\Gamma_c(Q_+)$ of
$(Q_+,d\theta\wedge \frac{dx}{x},\frac{dx}{x})$ applying Lemma
\ref{lcs} (choosing $\tilde{g}=\log x$, so that
$e^{-\tilde{g}}\tilde{\Omega}=d(x^{-1}d\theta)$ there). We obtain
the quotient of
$$(\tilde{Q}_+\times \RR \times \tilde{Q}_+,x_2
d\epsilon-\frac{x_2}{x_1}d\theta_1+d\theta_2,\frac{x_2}{x_1})$$
by the diagonal $\ZZ$ action on the variables
$(\theta_1,\theta_2)$. Here $(\theta_i,x_i)$ are the coordinates
on the two copies of the universal cover $\tilde{Q}_+\cong
\RR\times\RR_+$ and $\epsilon$ is the coordinate on the $\RR$
factor. The groupoid structure is given by the product of the pair
groupoid over $\tilde{Q}_+$ and group $\RR$. This contact
groupoid, and the one belonging to $Q_-=S^1\times \RR_-$, will sit
as open contact subgroupoids in the contact groupoid of $Q$, and
the question is how to ``complete'' the disjoint union of
${\Gamma}_c(Q_+)$ and ${\Gamma}_c(Q_-)$  to obtain the contact
groupoid of $Q$. A clue comes from the simplest case of groupoid
with two open orbits and a closed one to separate them, namely the
transformation groupoid of a vector field on $\RR$ with exactly
one zero. The transformation groupoid associated to $-x\partial_x$
is $\RR\times \RR=\{(t,x)\}$ with source given by $x$, target
given by $e^{-t}x$ and multiplication
$(t',e^{-t}x)\cdot(t,x)=(t'+t,x)$. Notice that, on each of the two
open orbits $\RR_+$ and $\RR_-$ the groupoid is isomorphic to a
pair groupoid by the correspondence $(t,x)\in \RR \times\RR_{\pm}
\mapsto (e^{-t}x,x)\in \RR_{\pm}\times \RR_{\pm}$, with inverse
$(x_1,x_2) \mapsto (\log(\frac{x_2}{x_1}),x_2)$.

Now we embed ${\Gamma}_c(Q_+)$ into the groupoid ${\Gamma}_c(Q)$
described in \eqref{groidgstar} by the mapping
$$(\theta_1,x_1,\epsilon,\theta_2,x_2)\mapsto
\left(\theta_1,t=\log(\frac{x_2}{x_1}),\epsilon,\theta_2,x=x_2\right),$$
and similarly for $\Gamma_c(Q_-)$. The contact forms and function
translate to those indicated in \eqref{groidgstar}, which as a
consequence also satisfy the multiplicativity condition. One
 checks directly that the one form is a contact form also on the complement
 $\{x=0\}$ of the two open subgroupoids.
Therefore the one described in \eqref{groidgstar} is a contact
groupoid, and since we know that the source map is a Jacobi map on
the open dense set sitting over $Q_+$ and $Q_-$, it is the contact
groupoid of $(Q,E\wedge x
\partial_x,E)$.

Now we consider the $S^1$ contact reduction of the above s.s.c.
groupoid $\Gamma_c(Q)$.  As shown in the proof of Theorem
\ref{red} the moment map is $J_{\Gamma}=1-f_{\Gamma}=1-e^t$, so its zero
level set is $\{t=0\}$. The definition of moment map and the fact
that the infinitesimal generator $v_{\Gamma}$ of the $S^1$ action
projects to $E$ both via source and via target imply that on
$\{t=0\}$ we have
$v_{\Gamma}=(\partial_{\theta_1},0,0,\partial_{\theta_2},0)$.
 So $J^{-1}(0)/S^1$ is $\RR^3$ with
coordinates $(\theta:=\theta_2-\theta_1, \epsilon,x)$, 1-form
$d\theta+xd\epsilon$, source and target both given by $x$ and
groupoid multiplication given by addition in the $\theta$ and
$\epsilon$ factors. Upon division of the $\theta$ factor by $\ZZ$
(notice that the Reeb vector field of $\Gamma_c(Q)$ is
$\partial_{\theta_2}$)
 this is clearly just the
prequantization of $T^*\RR$, endowed with the canonical symplectic
form $dx\wedge d\epsilon$ and fiber addition as groupoid
multiplication, i.e. the prequantization of the symplectic
groupoid of the Poisson manifold $(\RR,0)$.
\end{ep}

\appendix
\section{Lie algebroids of precontact groupoids}\label{preciso}

\begin{lemma}\label{idef}
Let $(\Gamma,\theta_{\Gamma},f_{\Gamma})$ be a precontact
groupoid (as in Definition \ref{prect}) over the Jacobi-Dirac
manifold $(Q,\bL)$, so that the source map be a Jacobi-Dirac map.
Then a Lie algebroid isomorphism between $\ker \bs_*|_Q$ and $\bL$
is given by
\begin{eqnarray}\label{sixteen}
Y \mapsto
(\bt_*Y,-{r_{\Gamma}}_*Y)\oplus(-d\theta_{\Gamma}(Y)|_{TQ},
\theta_{\Gamma}(Y))
\end{eqnarray}
where $e^{-r_{\Gamma}}=f_{\Gamma}$. A Lie algebroid isomorphism
between $\ker \bt_*|_Q$ and $\bL$ (obtained composing the above
with $i_*$ for $i$ the inversion) is
\begin{eqnarray} \label{seventeen}Y \mapsto
(\bs_*Y,{r_{\Gamma}}_*Y)\oplus(d\theta_{\Gamma}(Y)|_{TQ},
-\theta_{\Gamma}(Y))
\end{eqnarray}

\end{lemma}

\begin{proof}
Consider the groupoid $\Gamma\times\RR$ over $Q\times \RR$ with
target map $\tilde{\bt}(g,t)=(\bt(g),t-r_{\Gamma}(g))$ and the
obvious source $\tilde{\bs}$ and multiplication. $(\Gamma \times
\RR, d(e^t \theta_{\Gamma}))$ is then a presymplectic groupoid
with the property that $\tilde{\bs}$ is a forward Dirac map onto
$(Q\times \RR,\tilde{L})$, where
$$\tilde{L}_{(q,t)}=\{(X,f)\oplus e^t(\xi,g): (X,f)\oplus (\xi,g)\in L_q\}$$
 is the ``Diracization'' (\cite{WZ}\cite{IW}) of the Jacobi-Dirac
structure $\bar{L}$ and $t$ is the coordinate on $\RR$. In the
special case that $\bL$ corresponds to a Jacobi structure this is
just Prop. 2.7 of \cite{cz}; in the general case (but assuming
different conventions for the multiplicativity of
$\theta_{\Gamma}$ and for which of source and target is a
Jacobi-Dirac map) this is Prop. 3.3 in \cite{IW}. We will prove
only the first isomorphism above (the one for $\ker \bs_*|_Q$);
the other one follows by composing the first isomorphism with
$i_*$.
 Now we consider the
following diagram of spaces of sections (on the left column we
have sections over $Q$, on the right column sections over $Q\times
\RR$):
\[
\begin{CD}
\Gamma({\ker \bs_*}|_{Q}) @>\Phi_{\bs}>> \Gamma({\ker \tilde{\bs}_*}|_{Q\times \RR})\\
 @VVV \Phi @VVV\\
\Gamma(\bar{L}) @>\Phi_L>> \tilde{L}.
\end{CD}
\]
The first horizontal arrow $\Phi_{\bs}$ is $Y \mapsto \tilde{Y}$,
where the latter denotes the constant extension of $Y$ along the
$\RR$ direction of the base $Q\times \RR$.
 Notice that the projection $pr:\Gamma\times\RR \rightarrow \Gamma$
is a groupoid morphism, so it induces a surjective Lie algebroid
morphism $pr_*:\ker {\tilde{\bs}_*}|_{Q\times \RR} \rightarrow
\ker{\bs_*}|_Q$. Since sections $\tilde{Y}$ as above are
projectable, by Prop. 4.3.8. in \cite{MK2} we have
$pr_*[\tilde{Y}_1,\tilde{Y}_2]=[Y_1,Y_2]$, and since $pr_*$ is a
fiberwise isomorphism
% so that the right-invariant extensions of
%$Y$ on $\Gamma\times\RR$ (which is clearly projectable) is mapped
%into the right-invariant extension of $Y$. Since the Lie bracket
%of projectable vector fields is the Lie bracket of the projection
%\mcomment{is this reallly true?},
we deduce that $\Phi_{\bs}$ is a bracket-preserving map.

 The vertical arrow $\Phi$ is induced from the following
isomorphism of Lie algebroids (Cor. 4.8 iii of
\cite{bcwz}\footnote{In \cite{bcwz} the authors adopted the
convention that the target map be a Dirac map. Here we use their
result applied to the pre-symplectic form $-\Omega$.}
% and composing
%with the Lie algebroid isomorphism $\ker \tilde{\bt}_*\cong \ker
%\tilde{\bs}_*$ given by $i_*$, where $i$ is the inversion map, we
%obtain the desired formula.}
) valid for any presymplectic manifold $(\tilde{\Gamma},\Omega)$
over a Dirac manifold $(N,\tilde{L})$ for which the source map is
Dirac:
$$\ker {\tilde{\bs}_*}|_{N}\rightarrow \tilde{L}\;,\;
Z\mapsto(\tilde{\bt}_*Z, -\Omega(Z)|_{TN}).$$ In our case, as
mentioned above, the presymplectic form is $ d(e^t
\theta_{\Gamma})$.

The second horizontal arrow $\Phi_L$ is the natural map
$$(X,f)\oplus (\xi,g)\in L_q\mapsto (X,f)\oplus e^t(\xi,g)\in \bar{L}_{(q,t)}$$
 which
preserves the Lie algebroid bracket (see the remarks after
Definition 3.2 of \cite{WZ}).

One can check that $(\Phi
\circ\Phi_{\bs})(Y)=(\tilde{\bt}_*\tilde{Y})\oplus(-d(e^t
\theta_{\Gamma})(\tilde{Y})|_{TQ\times\RR})$
 lies in the
image of the injective map $\Phi_L$. The resulting map from
$\Gamma(\ker \bs_*)$ to $\Gamma(\bar{L})$ is given by
\eqref{sixteen} and the arguments above show that this map
preserves brackets. Further it is clear that this map of sections
is induced by a vector bundle morphism given by the same formula,
which clearly preserves not only the bracket of sections but also
the anchor, so that the map $\ker {\bs_*}|_Q\rightarrow \bL$ given
by \eqref{sixteen}
 is a Lie algebroid morphism.

To show that it is an isomorphism one can argue noticing that
 $\ker \bs_*$ and $\bL$ have the same dimension and
%and that the induced
 %map on sections \label{secsmap} is injective. Alternatively, one can
 show that the vector bundle map \label{secsmap} is injective,
by using the ``non-degeneracy condition'' in Def. \ref{prect}
 and the fact that the source and target fibers
 of $\Gamma\times\RR$ are pre-symplectic orthogonal to each other.
\end{proof}

The vector bundle morphisms in the above lemma give a
characterization of vectors tangent to the $\bs$ or $\bt$ fibers
of a precontact groupoid as follows.
 Consider
 for instance a vector $\lambda$ in $\bar{L}_x$, where $\bar{L}$ is the
 Jacobi-Dirac structure on the base $Q$. This vector
 corresponds to some
 $Y_x\in \ker\bt_*$  by the isomorphism \eqref{seventeen}, and by
 left translation we obtain
 a vector field $Y$ tangent to $\bt^{-1}(x)$.
Of course, every vector tangent to $\bt^{-1}(x)$ arises in this
way for a unique $\lambda$. The vector field $Y$ satisfies the
following equations at every point $g$ of $\bt^{-1}(x)$, which
follow by simple computation from the multiplicativity of
$\theta_{\Gamma}$:
 $\theta_{\Gamma}(Y_g)=\theta_{\Gamma}(Y_x)$,
$d\theta_{\Gamma}(Y_g,Z)=d\theta_{\Gamma}(Y_x,\bs_*Z)-{r_{\Gamma}}
_*Y_x\cdot\theta_{\Gamma}(Z)$ for all $Z\in T_g\Gamma$,
${r_{\Gamma}}_*Y_g={r_{\Gamma}}_*Y_x$ and $\bs_*Y_g=\bs_*Y_x$.
Notice that the right hand sides of this properties can be
expressed in terms of the four components of $\lambda \in
\cE^1(Q)$, and that by the ``non-degeneracy'' of $\theta_{\Gamma}$
these properties are enough to uniquely determine $Y_g$. We sum up
this discussion into the following corollary, which can be used as
a tool in computations on precontact groupoids in the same way
that hamiltonian vector fields are used on contact or symplectic
groupoids (such as the proof of Thm. \ref{red}):

\begin{cor}\label{computation}
Let $(\Gamma,\theta_{\Gamma},f_{\Gamma})$ be a precontact
groupoid (as in Definition \ref{prect}) and denote by $\bar{L}$
 the Jacobi-Dirac structure on the base $Q$ so that source map is
 Jacobi-Dirac. Then there is bijection between sections of $\bar{L}$ and vector
fields on $\Gamma$ which are tangent to the $\bt$-fibers and are
left invariant. To a section $(X,f)\oplus(\xi,g)$ of
$\bar{L}\subset \cE^1(Q)$ corresponds the unique vector field $Y$
tangent to the $\bt$-fibers which satisfies
\begin{itemize}
\item $\theta_{\Gamma}(Y)=-g$\\
\item $d\theta_{\Gamma}(Y)=\bs^*\xi-f\theta_{\Gamma}$\\
\item $\bs_*Y=X$.
\end{itemize}
$Y$ furthermore satisfies ${r_{\Gamma}}_*Y=f$.
\end{cor}

\section{Groupoids of locally conformal symplectic  manifolds }\label{lcsgroid}

A locally conformal symplectic (l.c.s.) manifold is a manifold
$(Q,\Omega,\omega)$ where $\Omega$ is a non-degenerate 2-form and
$\omega$ is a closed 1-form satisfying $d\Omega =\omega\wedge
\Omega$. Any Jacobi manifold is foliated by contact and l.c.s.
leaves (see for example \cite{zz}); in particular a l.c.s.
manifold is a Jacobi manifold, and hence, when it is integrable,
it has an associated s.s.c. contact groupoid. In this appendix we
will construct explicitly this groupoid; we make use of it in
Example \ref{1dim}.

\begin{lemma}\label{lcs}
Let $(Q,\Omega,\omega)$ a locally conformal symplectic manifold.
Consider the pullback structure on the universal cover
$(\tilde{Q},\tilde{\Omega},\tilde{\omega})$, and write
$\tilde{\omega}=d\tilde{g}$. Then $Q$ is integrable as a Jacobi
manifold iff the symplectic form $e^{-\tilde{g}}\tilde{\Omega}$ is
a multiple of an integer form. In that case, choosing $\tilde{g}$
so that $e^{-\tilde{g}}\tilde{\Omega}$ is integer,
 the s.s.c. contact groupoid of $(Q,\Omega,\omega)$ is
 the quotient of
\begin{eqnarray}\label{cgtq}\left(\tilde{R}\times_{\RR}\tilde{R},e^{\tilde{\bs}^*\tilde{g}}(-\tilde{\sigma}_1
+\tilde{\sigma}_2),\frac{e^{\tilde{\bs}^*\tilde{g}}}{e^{\tilde{\bt}^*\tilde{g}}}\right),
\end{eqnarray}
a groupoid over $\tilde{Q}$, by a natural $\pi_1(Q)$ action. Here
$(\tilde{R},\tilde{\sigma})$ is the universal cover (with the
pullback 1-form) of  a prequantization $(R,\sigma)$ of
$(\tilde{Q},e^{-\tilde{g}}\tilde{\Omega})$, and the group $\RR$
acts by the diagonal lift of the $S^1$ action on $R$.
\end{lemma}
\begin{proof}
% It is easy to see that the
%contact groupoid of $(Q,e^{-f}\Omega)$ is
%$(R\times_{S^1}R,\bs^*\sigma -\bt^*\sigma,1)$.
Using for example the Lie algebroid integrability criteria of
\cite{cf}, one sees that $(Q,\Omega,\omega)$ is integrable as a
Jacobi manifold iff $(\tilde{Q},\tilde{\Omega},\tilde{\omega})$
is.
 Lemma 1.5 in
Appendix I of \cite{zz} states that, given a contact groupoid,
multiplying the contact form by $\bs^*u$ and the multiplicative
function by $\frac{\bs^*u}{\bt^*u}$ gives another contact
groupoid, for any non-vanishing function $u$ on the base. Such an
operation corresponds to twisting the groupoid, viewed just as a
Jacobi manifold, by the function $\bs^*{u^{-1}}$, hence the Jacobi
structure induced on the base by the requirement that the source
be a Jacobi map is the twist of the original one by $u^{-1}$. So
$(\tilde{Q},\tilde{\Omega},\tilde{\omega})$  is integrable iff the
symplectic manifold $(\tilde{Q},e^{-\tilde{g}}\tilde{\Omega})$ is
Jacobi integrable, and by Section 7 of \cite{cz} this happens
exactly when the class of $e^{-\tilde{g}}\tilde{\Omega}$ is a
multiple of an integer one.

Choose $\tilde{g}$ so that this class is actually integer. A
contact groupoid of $(\tilde{Q},e^{-\tilde{g}}\tilde{\Omega})$ is
clearly $(R\times_{S^1}R,[-\sigma_1+\sigma_2],1)$, where the $S^1$
action on $R\times R$ is diagonal and ``$[\;\;]$'' denotes the
form descending from $R\times R$. This groupoid is not s.s.c.; the
s.s.c. one is $\tilde{R}\times_{\RR}\tilde{R}$, where the $\RR$
action on $\tilde{R}$ is the lift of the $S^1$ action on $R$. The
source simply connectedness follows since $\RR$ acts transitively
(even though not necessarily freely) on each fiber of the map
$\tilde{R}\rightarrow \tilde{Q}$, and this in turns holds because
any $S^1$ orbit in $R$ generates $\pi_1(R)$ and because the
fundamental group of a space always acts (by lifting loops)
transitively on the fibers of its universal cover.

By the above cited Lemma from \cite{zz} we conclude that
\eqref{cgtq} is the s.s.c. contact groupoid of
$(\tilde{Q},\tilde{\Omega},\tilde{\omega})$. The fundamental group
of $Q$ acts on $\tilde{Q}$ respecting its geometric structure, so
it acts on its Lie algebroid $T^*\tilde{Q}\times \RR$. Since the
path-space construction of the s.s.c. groupoid is canonical (see
Subsection \ref{dirid}), $\pi_1(Q)$ acts on the s.s.c. groupoid
\eqref{cgtq} preserving the groupoid and geometric structure.
Hence the quotient is a s.s.c. contact groupoid over
$(Q,\Omega,\omega)$, and its source map is a Jacobi map, so it is
the s.s.c. contact groupoid of $(Q,\Omega,\omega)$.
\end{proof}

\section{On a construction of Vorobjev}\label{vor}

In Section \ref{cpreq} we derived the
geometric structure on the circle bundles $Q$ from a
prequantizable Dirac manifold $(P,L)$ and a suitable choice of
connection $D$. In this appendix we describe an alternative attempt;
even though we can make our construction work only
if we start with a symplectic manifold, we believe the
construction is interesting on its own right.

First we recall Vorobjev's construction in Section 4 of \cite{Vo},
which the author there uses to study the linearization problem of
Poisson manifolds near a symplectic leaf. Consider a transitive
algebroid $A$ over a base $P$ with anchor $\rho$;  the kernel
$\ker \rho$ is a bundle of Lie algebras.
 Choose a splitting $\gamma:TP\rightarrow A$ of the anchor.
 Its curvature $R_{\gamma}$ is a 2-form on $P$ with values in
 $\Gamma(\ker \rho)$  (given by $R_{\gamma}(v,w)=[\gamma v, \gamma w]_A-
\gamma[v,w])$. The splitting $\gamma$ also induces a
(TP-)covariant derivative $\nabla$ on $\ker \rho$ by
$\nabla_vs=[\gamma v,s]_A$ . Now, if $P$ is endowed with a
symplectic form $\omega$, a neighborhood of the zero section in
$(\ker \rho)^*$ inherits a Poisson structure
$\Lambda_{vert}+\Lambda_{hor}$ as follows (Theorem 4.1 in
\cite{Vo}): denoting by $F_s$ the fiberwise linear function on
$(\ker \rho)^*$ obtained by contraction with the section $s$ of
$\ker \rho$, the Poisson bivector has a vertical component
determined by $\Lambda_{vert}(dF_{s_1},dF_{s_2})=F_{[s_1,s_2]}$.
It also has a component  $\Lambda_{hor}$  which is tangent to the
Ehresmann connection $Hor$ given by the dual
connection\footnote{In \cite{Vo} the author phrases this condition
as $\cL_{hor(X)}F_s=F_{\nabla_Xs}$.} to $\nabla$ on the bundle
$(\ker \rho)^*$; $\Lambda_{hor}$ at $e\in (\ker \rho)^*$ is
obtained by restricting the non-degenerate form $\omega-\langle
R_{\gamma}, e \rangle $ to $Hor_e$  and inverting it. (Here we are
identifying $Hor_e$ and the corresponding tangent space to $P$.)

To apply Vorobjev's construction in our setting, let $(P,\omega)$
be a prequantizatible symplectic manifold and $(K,\nabla_K)$ its
prequantization line bundle with Hermitian connection of curvature
$2\pi i \omega$. By Lemma \ref{flat} we obtain a flat
$TP\oplus_{\omega}\RR$-connection $\tilde{D}_{(X,f)}=\nabla_X+2\pi
i f$ on $K$. Now we make use of the following well know fact about
extensions, which can be proven by direct computation:
\begin{lemma}
Let $A$ be a Lie algebroid over $M$, $V$ a vector bundle over $M$,
and $\tilde{D}$ a flat $A$-connection on $V$. Then $A\oplus V$
becomes a Lie algebroid with the anchor of $A$ as anchor and
bracket
$$[(Y_1,S_1),(Y_2,S_2)]=([Y_1,Y_2]_A,\tilde{D}_{Y_1}S_2-\tilde{D}_{Y_2}S_1).$$
\end{lemma}

Therefore $A:=TP\oplus_{\omega}\RR\oplus K$ is a transitive Lie
algebroid over $P$, with isotropy bundle $\ker \rho=\RR\oplus K$
and bracket $[(f_1,S_1),(f_2,S_2)]=[(0,2\pi i (f_1S_2-f_2S_1)]$
there.
 Now choosing the
canonical splitting $\gamma$ of the anchor
$TM\oplus_{\omega}\RR\oplus K \rightarrow TM$ we see that its
curvature is $R_{\gamma}(X_1,X_2)=(0,\omega(X_1,X_2),0)$. The
horizontal distribution on the dual of the isotropy bundle is the
product of the trivial one  on $\RR$ and of the one corresponding
to $\nabla_K$ on $K$ (upon identification of $K$ and $K^*$ by the
metric). By the above, there  is
 a Poisson
structure on $\RR\oplus K$, at least near the zero section: the
Poisson bivector at $(t,q)$ has a horizontal component given by
lifting the inverse of $(1-t)\omega$ and a vertical component
which turns out to be $2\pi (i q\partial_q)\wedge
\partial t$, where ``$iq\partial_q$'' denotes the vector field tangent to the
circle bundles in $K$ obtained by turning by $90^{\circ}$ the
Euler vector field $q\partial_q$. A symplectic leaf is clearly
given by $\{t<1\}\times Q$ (where $Q=\{|q|=1\}$). On this leaf the
symplectic structure is seen to be given by
$(1-t)\omega+\theta\wedge dt=d((1-t)\theta)$, where $\theta$ is
the connection 1-form on $Q$ corresponding to the connection
$\nabla_K$ on $K$ (which by definition satisfies
$d\theta=\pi^*\omega$). This means that the leaf is just the
symplectification $(\RR_+\times Q,d(r\theta))$
 of $(Q,\theta)$ (here $r=1-t$), which is a ``prequantization
 space'' for our symplectic manifold $(P,\omega)$.
Unfortunately we are not able to modify Vorobjev's construction
appropriately when $P$ is a Poisson or Dirac manifold.

\bibliographystyle{habbrv}
\bibliography{bibz}

\def\cprime{$'$} \def\cprime{$'$}
\begin{thebibliography}{10}

\bibitem{albert}
C.~Albert.
\newblock Le th\'eor\`eme de r\'eduction de {M}arsden-{W}einstein en
  g\'eom\'etrie cosymplectique et de contact.
\newblock {\em J. Geom. Phys.}, 6(4):627--649, 1989.

\bibitem{bcz}
F.~Bonechi, A.~S. Cattaneo, and M.~Zabzine.
\newblock {Geometric quantization and non-perturbative Poisson sigma model},
  arxiv:math.SG/0507223.

\bibitem{bcwz}
H.~Bursztyn, M.~Crainic, A.~Weinstein, and C.~Zhu.
\newblock Integration of twisted {D}irac brackets.
\newblock {\em Duke Math. J.}, 123(3):549--607, 2004.

\bibitem{cmdl}
D.~Chinea, J.~C. Marrero, and M.~de~Le{\'o}n.
\newblock Prequantizable {P}oisson manifolds and {J}acobi structures.
\newblock {\em J. Phys. A}, 29(19):6313--6324, 1996.

\bibitem{Co}
T.~J. Courant.
\newblock Dirac manifolds.
\newblock {\em Trans. Amer. Math. Soc.}, 319(2):631--661, 1990.

\bibitem{Cr}
M.~Crainic.
\newblock {Prequantization and Lie brackets}, arxiv:math.DG/0403269.

\bibitem{cf}
M.~Crainic and R.~L. Fernandes.
\newblock Integrability of {L}ie brackets.
\newblock {\em Ann. of Math. (2)}, 157(2):575--620, 2003.

\bibitem{cz}
M.~Crainic and C.~Zhu.
\newblock {Integrability of Jacobi structures, math.DG/0403268, to appear in
  Annal of Fourier Institute}.

\bibitem{tudor}
O.~Dr{\u{a}}gulete, L.~Ornea, and T.~S. Ratiu.
\newblock Cosphere bundle reduction in contact geometry.
\newblock {\em J. Symplectic Geom.}, 1(4):695--714, 2003.

\bibitem{for}
R.~Fernandes, J.~Ortega, and T.~Ratiu.
\newblock {Momentum maps in Poisson geometry}.
\newblock in preparation.

\bibitem{Fe}
R.~L. Fernandes.
\newblock Lie algebroids, holonomy and characteristic classes.
\newblock {\em Adv. Math.}, 170(1):119--179, 2002.

\bibitem{Hu}
J.~Huebschmann.
\newblock Poisson cohomology and quantization.
\newblock {\em J. Reine Angew. Math.}, 408:57--113, 1990.

\bibitem{IMa}
D.~Iglesias and J.~Marrero.
\newblock {Lie algebroid foliations and $E^1(M)$-Dirac structures},
  arXiv:math.DG/0106086.

\bibitem{kost}
B.~Kostant.
\newblock Quantization and unitary representations. {I}. {P}requantization.
\newblock In {\em Lectures in modern analysis and applications, III}, pages
  87--208. Lecture Notes in Math., Vol. 170. Springer, Berlin, 1970.

\bibitem{MK2}
K.~C.~H. Mackenzie.
\newblock {\em General theory of {L}ie groupoids and {L}ie algebroids}, volume
  213 of {\em London Mathematical Society Lecture Note Series}.
\newblock Cambridge University Press, Cambridge, 2005.

\bibitem{MW}
K.~Mikami and A.~Weinstein.
\newblock Moments and reduction for symplectic groupoids.
\newblock {\em Publ. Res. Inst. Math. Sci.}, 24(1):121--140, 1988.

\bibitem{IW}
D.~I. Ponte and A.~Wade.
\newblock {Integration of Dirac-Jacobi structures}, arXiv:math.DG/0507538.

\bibitem{So}
J.-M. Souriau.
\newblock Quantification g\'eom\'etrique.
\newblock In {\em Physique quantique et g\'eom\'etrie (Paris, 1986)}, volume~32
  of {\em Travaux en Cours}, pages 141--193. Hermann, Paris, 1988.

\bibitem{tz}
H.-H. Tseng and C.~Zhu.
\newblock Integrating {L}ie algebroids via stacks.
\newblock {\em Compos. Math.}, 142(1):251--270, 2006.

\bibitem{va}
I.~Vaisman.
\newblock On the geometric quantization of {P}oisson manifolds.
\newblock {\em J. Math. Phys.}, 32(12):3339--3345, 1991.

\bibitem{Vo}
Y.~Vorobjev.
\newblock {Coupling Tensors and Poisson Geometry Near a Single Symplectic
  Leaf}, arxiv:math.SG/0008162.

\bibitem{Wa}
A.~Wade.
\newblock Conformal {D}irac structures.
\newblock {\em Lett. Math. Phys.}, 53(4):331--348, 2000.

\bibitem{w-ncomquan}
A.~Weinstein.
\newblock Noncommutative geometry and geometric quantization.
\newblock In {\em Symplectic geometry and mathematical physics
  (Aix-en-Provence, 1990)}, volume~99 of {\em Progr. Math.}, pages 446--461.
  Birkh\"auser Boston, Boston, MA, 1991.

\bibitem{wx}
A.~Weinstein and P.~Xu.
\newblock Extensions of symplectic groupoids and quantization.
\newblock {\em J. Reine Angew. Math.}, 417:159--189, 1991.

\bibitem{WZ}
A.~Weinstein and M.~Zambon.
\newblock Variations on prequantization.
\newblock In {\em Travaux math\'ematiques. Fasc. XVI}, Trav. Math., XVI, pages
  187--219. Univ. Luxemb., Luxembourg, 2005.

\bibitem{willett}
C.~Willett.
\newblock Contact reduction.
\newblock {\em Trans. Amer. Math. Soc.}, 354(10):4245--4260 (electronic), 2002.

\bibitem{zz}
M.~Zambon and C.~Zhu.
\newblock Contact reduction and groupoid actions.
\newblock {\em Trans. Amer. Math. Soc.}, 358(3):1365--1401 (electronic), 2006.

\bibitem{z:lie2}
C.~Zhu.
\newblock {Lie II theorem for Lie algebroids via stacky Lie groupoids},
  arXiv:math.DG/0701024.

\end{thebibliography}

\noindent Chenchang Zhu\\
 Department Mathematik, Eidgen\"ossische Technische Hochschule (ETH)\\
 R\"amistr. 101, 8092 Z\"urich, Switzerland\\
 zhu@math.ethz.ch

\vspace*{0.8cm}

\noindent Marco Zambon\\
      Mathematisches Institut, Universit\"at Z\"urich\\
      Winterthurerstr. 190, 8057 Z\"urich, Switzerland\\
  zambon@math.unizh.ch

\end{document}